\documentclass[11]{article}
\usepackage{amsmath,amssymb,amsbsy}
\usepackage{xspace}
\usepackage{url}
\usepackage[inline]{enumitem}

\usepackage{xifthen}
\usepackage{proof}

\setlist{leftmargin=*}
\setitemize[2]{label=$\bullet$}
\setenumerate[1]{label=(\arabic*), ref=(\arabic*)}
\setenumerate[2]{leftmargin=*,label=(\alph*),ref=(\theenumi.\alph*)}
\setenumerate[3]{label=(\roman*)), ref=(\theenumii.\roman*)}

\newcommand{\F}{{\bf F}}
\newcommand{\K}{{\bf K\,}\xspace}

\newcommand{\rmR}{\mathrm{R}}

\newcommand{\calL}{\mathcal{L}}

\newcommand{\calS}{\mbox{$\mathcal{S}$}\xspace}

\newcommand{\calT}{\mathcal{T}}
\newcommand{\calC}{\mathcal{C}\xspace}

\newcommand{\calM}{\mathcal{M}\xspace}

\newcommand{\calR}{\mathcal{R}\xspace}

\renewcommand{\a}{\alpha}
\renewcommand{\b}{\beta}
\newcommand{\g}{\gamma}

\newcommand{\s}{\sigma}

\newcommand{\D}{\Delta}
\newcommand{\G}{\Gamma}
\newcommand{\T}{\mbox{$\Theta$}}

\newcommand{\forza}{\Vdash}
\newcommand{\sat}{\triangleright} 
\newcommand{\nonforza}{\nVdash}

\newtheorem{theorem}{Theorem}

\newtheorem{proposition}{Proposition}
\newtheorem{example}{Example}

\def\squareforqed{\hbox{\rlap{$\sqcap$}$\sqcup$}}
\def\qed{\ifmmode\squareforqed\else{\unskip\nobreak\hfil
\penalty50\hskip1em\null\nobreak\hfil\squareforqed
\parfillskip=0pt\finalhyphendemerits=0\endgraf}\fi}

\newcommand{\uK}{{\cal K}}
\newcommand{\uM}{{\cal M}}

\newcommand{\kripke}{\mbox{$\langle S$,} \mbox{$\rho$,} \mbox{$\leq,$} \mbox{$\E$,}
  \mbox{$\forza\rangle$}\xspace}

\newcommand{\PV}{\mathrm{V}}

\newcommand{\at}{\mathrm{At}}
\newcommand{\AT}{\at}
\newcommand{\depth}{\mathrm{depth}}
\newcommand{\md}{\mathrm{maxdepth}}
\newcommand{\noninv}{\mathrm{nonInv}}
\newcommand{\inv}{\mathrm{Inv}}
\newcommand{\mindepth}{\mathrm{mindepth}}

\def\regola#1#2#3{\infer[(#3)]{#1}{#2}}
\def\apply#1#2#3{\textrm{$\displaystyle\frac{#1}{#2}\,\,\textstyle #3$}}

\newcommand{\seq}[3]{
  \ifthenelse{\isempty{#1}}{\emptyset}{#1}\, ; \,
  \ifthenelse{\isempty{#2}}{\emptyset}{#2}\, \Longrightarrow
  \ifthenelse{\isempty{#3}}{\bot}{#3}
  }

  \newcommand{\eseq}[3]{
    \ifthenelse{\isempty{#1}}{\emptyset}{#1}\, ; \,
    \ifthenelse{\isempty{#2}}{\emptyset}{#2}\,
    \Longrightarrow
    \kern-3.5ex\raise 1.3ex\hbox{\sc e}\kern 2.7ex
    \ifthenelse{\isempty{#3}}{\bot}{#3}
  }

  \newcommand{\mseq}[3]{
    \ifthenelse{\isempty{#1}}{\emptyset}{#1}\, ; \,
    \ifthenelse{\isempty{#2}}{\emptyset}{#2}\,
    \Longrightarrow
    \kern-3.5ex\raise 1.3ex\hbox{\sc m}\kern 2.7ex
    \ifthenelse{\isempty{#3}}{\bot}{#3}
  }

\newcommand{\false}{\bot}

\newlength{\brsep}
\newlength{\virsep}
\newlength{\pvsep} 
\newcommand{\iel}{{\bf IEL}\xspace}
\newcommand{\Int}{{\bf Int}\xspace}
\newcommand{\ielmeno}{{\bf IEL$^-$}\xspace}

\newcommand{\iea}{{$\K(A\to B) \to (\K A\to \K B)$}}
\newcommand{\ieb}{$A\to\K A$}
\newcommand{\iec}{$\K A\to\neg\neg A$}

\newcommand{\E}{\textrm{$\,\mathrm E\,$}\xspace}

\newcommand{\LSJ}{\mbox{$\mathrm{LSJ}$}\xspace}
\newcommand{\Liel}{\mbox{$\mathrm{Liel}$}\xspace}
\newcommand{\proc}{\textrm{$\mathrm{Piel}$}\xspace}
\newcommand{\priel}{\textrm{$\mathrm{PR}$}\xspace}
\newcommand{\Riel}{\mbox{$\mathrm{Riel}$}\xspace}

\DeclareMathAlphabet\mathbfcal{OMS}{cmsy}{b}{n}

\begin{document}

\title{Linear Depth Deduction with Subformula Property for Intuitionistic Epistemic Logic}

\author{
  Guido Fiorino\\
  DISCO, Universit\`a degli Studi di Milano-Bicocca,\\
  Viale Sarca, 336, 20126, Milano, Italy.\\{\tt guido.fiorino@unimib.it}}


\maketitle

\begin{abstract}
  In their seminal paper Artemov and Protopopescu provide Hilbert formal systems,
  Brower-Heyting-Kolmogorov and Kripke semantics for the logics of intuitionistic belief
  and knowledge.  Subsequently Krupski has proved that the logic of intuitionistic
  knowledge is PSPACE-complete and Su and Sano have provided calculi enjoying the
  subformula property. This paper continues the investigations around to sequent calculi
  for Intuitionistic Epistemic Logics by providing sequent calculi that have the
  subformula property and that are terminating in linear depth. Our calculi allow us to
  design a procedure that for invalid formulas returns a Kripke model of minimal
  depth. Finally we also discuss refutational sequent calculi, that is sequent calculi to
  prove the invalidity.
\end{abstract}

\makeatletter
\setlength\parskip{6\p@ \@plus \p@}
\setlength\parindent{0\p@}
\makeatother

\section{Introduction}
In~\cite{ArtPro:2016} the epistemic logics \iel and $\iel^-$ are introduced with the aim
to study the knowledge from the intuitionistic point of view. The authors remark that
Brower-Heyting-Kolmogorov (BHK) semantics is the intended semantics of intuitionistic logic, where a
proposition is true if it is proved and thus, in the
intuitionistic perspective, knowledge and belief are the product of verification.
From the idea that intuitionistic proof is a form of verification that implies
intuitionistic knowledge (represented by the modality $\K$) follows that the {\em
  co-reflexion principle} (or {\em constructivity of the proof}) $A\to \K A$ is assumed
both for belief and knowledge.
On the other hand, in the intuitionistic perspective, the verification of a statement does
not guarantee to a have a proof. A possibility is that a verified statement cannot be false
and this is the {\em intuitionistic reflexion principle} $\K A\to \neg \neg A$, formally
stating that if $A$ has a verification, then $A$ has a proof which is not necessarily
obtainable from the verification process. If the intuitionistic reflexion
principle is assumed, then known formulas cannot be false. If co-reflexion and
intuitionistic reflexion are assumed, then
intuitionistic truth implies intuitionistic knowledge and 
intuitionistic knowledge implies classical truth. The intuitionistic systems of epistemic
logic $\iel$ and $\iel^-$ differ because $\iel^-$ does not assume the
intuitionistic reflection principle thus we can have false beliefs.
The investigations in~\cite{ArtPro:2016} also include characterisations of \iel and
$\iel^-$ by means of both Hilbert axiom systems and Kripke semantics, where a binary
relation \E between worlds is added to the usual elements of the intuitionistic Kripke
frame.

Logics $\iel$ and
$\iel^-$ have attracted people engaged type theory: paper~\cite{PeriniBrogi:2020} provides
a formal analysis of the computational content of intuitionistic belief by introducing a
natural deduction calculus for $\iel^-$; paper~\cite{Rogozin:2020} constructs a type
system which is Curry-Howard isomorphic to $\iel^-$.
The investigations of our paper are more in line with those in~\cite{KruYat:2016}, where
it is proved that \iel is PSPACE-complete. Moreover, paper~\cite{KruYat:2016} extends Gentzen
calculus {\bf LJ} with two rules for the connective $\K$ that do not fulfil the
subformula property.
We also quote~\cite{SuSano:2019}, that presents a calculus for First Order \iel that
extends {\bf LJ} with new rules for the connective $\K$ and the result is a logical
apparatus that fulfil the subformula property.  Apart the aspect related to the
subformula property and computational complexity, in the quoted papers there is no
investigation about efficiency in proof search.

In this paper we propose sequent calculi
that have the subformula property and whose deductions are all depth-bounded in the number
of connectives occurring in the formula to be decided. To this aim, we use the ideas
from~\cite{FerFioFio:2013} where a sequent calculus for propositional intuitionistic logic
is provided. Paper~\cite{FerFioFio:2013} uses an extension of the ordinary sequent that
have similarities with nested sequents~\cite{Fitting:2014}.
Because of connective $\K$ and its Kripkean meaning, to handle \iel and $\iel^-$ we have
to extend the the object language used in~\cite{FerFioFio:2013} by adding a new type of
sequent that roughly speaking stores a semantical information related to the relation
\E. The sequent calculi we present can be explained by Kripke semantics of Intuitionistic
Epistemic Logics and our proofs follow model theoretic techniques, that is we provide
correctness and completeness theorems using the Kripke semantics for the logics at
hand.

The paper is organised as follows: in Section~\ref{sec:definitions} we recall the
definitions of \iel and $\iel^-$; in Section~\ref{sec:iel} we analyse the case of logic
\iel and we provide a sequent calculus, \Liel, and then, in
Section~\ref{sec:liel:completeness}, a procedure, \proc, that systematically builds trees
of sequents and returns a proof of linear depth, if the given formula is valid in \iel, or
a counter-model of minimal depth witnessing the invalidity of the given formula. In
Section~\ref{sec:riel}, we discuss \Riel, a calculus to prove that a formula is invalid in
\iel. Calculus \Riel is tied to \Liel: it has the subformula property and the deductions
are depth-bounded in the number of connectives occurring in the formula to be proved. In
Section~\ref{sec:lielmeno}, the sequent calculi $\Liel^-$ and $\Riel^-$ to prove validity
and invalidity of $\iel^-$ are provided. They are related to \Liel and $\Riel$
respectively and they are obtained by following the Kripke semantics of the logics. We
conclude with Section\ref{sec:conclusions} by discussing some possible future works.

\section{Definitions and Notations}
\label{sec:definitions}
Let $\PV$ be a denumerable set of {\em propositional variables}. We consider the propositional
language $\calL$ built  by using  the set of {\em atoms}
\mbox{$\at=\PV\cup\{\false\}$} and the set of connectives $\{\lor,\land,\to,\K\}$. When
convenient we write $\neg A$ in place of $A\to \bot$.

Logic \iel is proof-theoretically defined in~\cite{ArtPro:2016} as follows:
\begin{enumerate}[label=(Ax \arabic*)., ref=(Ax \arabic*)]
\item \label{Ax1}axioms for propositional intuitionistic logic $\Int$;
\item\label{Ax2} \iea;
\item \label{Ax3}\ieb;
\item \label{Ax4}\iec.
\end{enumerate}
In~\cite{ArtPro:2016} it is proved that \iel can be semantically characterised by (Kripke)
models \mbox{$\uK=\langle S,\rho,$} \mbox{$\leq,$} \mbox{$\E,$} \mbox{$\forza\rangle$}
defined as follows:
\begin{itemize}[label=-]
\item $\langle S,\rho,\leq\rangle$ is the usual Kripke frame
  for \Int;
\item $\E\subseteq S\times S$ fulfils the following properties:
  \begin{enumerate}[label=(Im \arabic*)., ref=(Im \arabic*)]
  \item\label{Im1}
    for every $\a,\b\in S,\ \a \E \b $ implies $\a\leq \b$;
  \item\label{Im2} for every $\a,\b,\g\in S$, if $\a\leq\b$ and $\b \E \g$, then $\a \E \g$;
  \item \label{E-prop}\label{Im3} for every $\a\in S$, there exists $\b\in S$ such that $\a \E \b$.
  \end{enumerate}

\item $\forza\subseteq S\times\calL$ is the {\em forcing relation} satisfying the
  following properties:
  \begin{itemize}[label=-]
  \item for every $p\in\PV$, for every $\a,\b\in S$, if $\a\forza p$ and $\a\leq \b$, then
    $\b\forza p$;
  \item for every $\a\in S$, $\a\nonforza \bot$;
  \item for every $\a\in S$, $\a\forza A\land B$ iff $\a\forza A$ and $\a\forza B$;
  \item for every $\a\in S$, $\a\forza A\lor B$ iff $\a\forza A$ or $\a\forza B$;
  \item for every $\a\in S$, $\a\forza A\to B$ iff for every $\b\in S$, if $\a\leq \b$,
    then $\b\nonforza A$ or $\b\forza B$;
  \item for every $\a\in S$, $\a\forza \K A$ iff for every $\b\in S$, if $\a\E \b$ then  $\b\forza A$.\
  \end{itemize}
\end{itemize}
Thus, $\iel=\{A\in\calL |$ for every Kripke model $\uK$, $\uK\forza A\}$.

By the properties of \E and $\forza$ it follows that the persistence property is fulfilled.

In the Kripke frames $\langle S, \rho,\leq\rangle$, the elements of $S$ are called {\em worlds},
or {\em states}, and in $S$ we distinguish the element $\rho$, the {\em root},
such that for every $\a\in S$, $\rho\leq \a$,
and the {\em final states}, where $\a\in S$ is a final state iff for every $\b\in S$, if $\a\leq \b$, then
$\a=\b$. Finally, given $\a,\b\in S$, we write $\a<\b$ to mean $\a\leq \b$ and $\a\neq \b$. 

To explain our ideas, we remark that by Property~\ref{E-prop} it follows that
every final state $\g$ of a model $\uK$ we have that $\g\E\g$ holds.
Thus, if $\a\forza \K A$, then necessarily all the final states reachable from
$\a$ force $A$ and if $\a\nonforza \K A$, then there exists at least a world $\b$ such that
$\a\E \b$ and $\b\nonforza A$ hold. By Property~\ref{Im1} follows that $\a\leq\b$ holds.
%
%
\newcommand{\ruleID}{(\mathbf{ Id})\xspace}
\newcommand{\ruleIRR}{(\mathbf{ Irr})\xspace}
\newcommand{\ruleeID}{(\mathbf{ eId})\xspace}
\newcommand{\ruleeIRR}{(\mathbf{ eIrr})\xspace}
\newcommand{\rulemIRR}{(\mathbf{ mIrr})\xspace}
\newcommand{\ruleSAT}{(\mathbf{ Sat})\xspace}
\newcommand{\ruleeSAT}{(\mathbf{ eSat})\xspace}

\newcommand{\AndLeft}{\mathbf{(\land L)}\xspace}
\newcommand{\AndRight}{\mathbf{(\land R)}\xspace}
\newcommand{\OrLeft}{\mathbf{(\lor L)}\xspace}
\newcommand{\OrRight}{\mathbf{(\lor R)}\xspace}
\newcommand{\ToLeft}{\mathbf{(\to L)}\xspace}
\newcommand{\ToRight}{\mathbf{\to R}\xspace}

\begin{figure}[t] 
  \[
    \begin{array}{c}
      {\bf Axioms} \\[4ex]
      \seq{\T}{\false,\G}{\D}~~\ruleIRR %
      \hspace{4em}
      \seq{\T}{A,\G}{A,\D}~~\ruleID %
      \\[4ex]
      {\bf Rules}\\[4ex]
      \infer[(\land L)]{\seq{\T}{A \land B,\G}{\D}}
      {\seq{\T}{A, B, \G}{\D}} %
      \hspace{4em}
      \infer[(\land R)]{\seq{\T}{\G}{A \land B, \D}}
      {\seq{\T}{\G}{A,\D} & \seq{\T}{\G}{B,\D}} %
      \\[4ex]
      \infer[(\lor L)]{\seq{\T}{A \lor B,\G}{\D}}
      {\seq{\T}{A,\G}{\D} & \seq{\T}{B,\G}{\D}} %
                            \hspace{4em}
                            \infer[(\lor R)]{\seq{\T}{\G}{A \lor B, \D}}
                            {\seq{\T}{\G}{A,B,\D}} %
      \\[4ex]
      \infer[(\to L)]{\seq{\T}{A \to B,\G}{\D}}
      {\seq{\T}{B,\G}{\D} & \seq{B,\T}{\G}{A,\D} & \seq{B}{\Theta,\G}{A}}
      \\[4ex]
      \regola{\seq{\T}{\G}{A \to B,\D}}
      {\seq{\T}{A, \G}{B,\D} & \seq{}{A,\T, \G}{B}}{\to R}
    \end{array}
  \]
  \caption{The calculus \LSJ for \Int.}
  \label{fig:lsj}
\end{figure}                     
\begin{figure}[ht]
  \small
  \[
    \setlength{\arraycolsep}{20pt}
    \begin{array}{c}
      \textrm{Rules and axioms in Figure~\ref{fig:lsj}}\medskip\\
      +
      \medskip\\
      {\bf Axioms} \\[4ex]
      \eseq{\T}{\false,\G}{\D}~~\ruleeIRR %
      \hspace{4em}
      \eseq{\T}{A,\G}{A,\D}~~\ruleeID %
      \\[4ex]
      \regola{\seq{\T}{\K A,\G}{\D}}
      {\eseq{\false}{A,\G}{\D}\ &\ \eseq{\false}{A, \T,\G}{\false}}{\K L} %
                                 \medskip\\
      \regola{\seq{\T}{\K A_1,\dots,\K A_n,\G }{\K B,\D}}
      {\eseq{\T}{A_1,\dots, A_n,\G}{B,\D}\ &\ \seq{}{A_1,\dots,A_n,\T,\G}{B}}{\K
                                            R}\smallskip\\
      \textrm{where $n\geq 0$}
      \medskip\\
      {
      \def\concl{\eseq{\T}{\K A, \G}{\D}}
      \def\prem{\eseq{\T}{A,\G}{\D}}
      \regola{\concl}{\prem}{e\K L}
      }
      \hspace*{\fill}
      {
      \def\concl{\eseq{\T}{\G}{\K B,\D}}
      \def\prem{\eseq{\T}{\G}{B,\D}}
      \regola{\concl}{\prem}{e\K R}
      }
      \medskip\\
      \infer[(e\land L)]{\eseq{\T}{A \land B,\G}{\D}}
      {\eseq{\T}{A, B, \G}{\D}} %
      \hspace{4em}
      \infer[(e\land R)]{\eseq{\T}{\G}{A \land B, \D}}
      {\eseq{\T}{\G}{A,\D} & \eseq{\T}{\G}{B,\D}} %
      \medskip\\
      \infer[(e\lor L)]{\eseq{\T}{A \lor B,\G}{\D}}
      {\eseq{\T}{A,\G}{\D} & \eseq{\T}{B,\G}{\D}} %
                            \hspace{4em}
                            \infer[(e\lor R)]{\eseq{\T}{\G}{A \lor B, \D}}
                            {\eseq{\T}{\G}{A,B,\D}} %
      \medskip\\
      \infer[(e\to L)]{\eseq{\T}{A \to B,\G}{\D}}
      {\eseq{\T}{B,\G}{\D} & \eseq{B,\T}{\G}{A,\D} & \seq{B}{\Theta,\G}{A}}
      \medskip\\
      \regola{\eseq{\T}{\G}{A \to B,\D}}
      {\eseq{\T}{A, \G}{B,\D} & \seq{}{A,\T, \G}{B}}{e\to R}
    \end{array}
  \]
  \caption{The calculus \Liel for \iel.}
  \label{fig:rules-for-k}
\end{figure}                     
\clearpage
Property~\ref{Im3} is crucial to the validity of~\ref{Ax4}. The logic $\iel^-$ is the
logic proof-theoretically characterised by the axioms~\ref{Ax1}-\ref{Ax3}. The
corresponding semantical characterisation is by Kripke models where relation \E
satisfies~\ref{Im1} and~\ref{Im2}.

The object language of our calculi are based on sequents of the kind $\seq{\T}{\G}{\D}$
and $\eseq{\T}{\G}{\D}$ with three compartments. We refer to the compartments respectively as {\em
  first, second} and {\em third compartment}.
For our purposes, the sets are always finite.
We call $\E$-sequents the sequents of
the kind $\eseq{\T}{\G}{\D}$.


The presence of the $\E$-sequents in the logical apparatus is related to the presence of
relation $\E$ in the definition of Kripke semantics for \iel. The proof of the
completeness theorem in the part related to rule $\K R$ uses the $\E$-sequents. Calculus
\LSJ of~\cite{FerFioFio:2013} handles sequents of the kind $\seq{\T}{\G}{\D}$ that are an
extension of the standard sequent $\G\Rightarrow\D$.  The meaning of the sequents can be
defined by means of Kripke models. Let $\uK=\kripke$ be Kripke model and $\a\in S$. We say
that {\em $\a$ satisfies ${\G}\Rightarrow{\D}$}, and we write $\a\sat {\G}\Rightarrow{\D}$, iff
for every $A\in\G$, $\a\forza A$ and for every $B\in\D$, $\alpha\nonforza B$.
We also say that {\em $\uK$ satisfies} $\G\Rightarrow\D$ and that {\em $\G\Rightarrow\D$ is satisfiable
(by $\uK$)}. 
For sake of
completeness, we recall that this corresponds to the following definition: $\a\sat\G\Rightarrow\D$ iff
$\a\forza \bigwedge\G\to\bigvee\D$, where $\bigwedge\G$ is the formula obtained by the
conjunction of all the formulas in $\G$ and $\bigvee\D$ is the formula obtained by the
disjunction of all the formulas in $\D$ (with the proviso that if $\G=\emptyset$ then
$\bigwedge \G$ is any tautological formula, such as $\top$, and if $\D=\emptyset$ then
$\bigvee \D$ is any contradictory formula, such as $\bot$).
We say that $\a$ satisfies $\seq{\T}{\G}{\D}$, and we write $\a\sat \seq{\T}{\G}{\D}$, iff
$\a\sat \G\Rightarrow\D$ and for every $C\in \T$ and for every $\b\in S$, if $\a<\b$, then
$\b\forza C$ (equivalently $\a\sat\seq{\T}{\G}{\F}$ iff $\a\sat \G\Rightarrow\D$ and for
every $\b\in S$, if $\a<\b$, then $\b\forza \bigwedge \T $).
The consequence of the semantical meaning of the sequents of the kind $\seq{\T}{\G}{\D}$
is that the rules handling the implication have one more premise than the rules handling
the standard sequent $\G\Rightarrow\D$, because the semantics of implication is defined
considering worlds $\b$ that are equal or greater than $\a$ and thus the rules take into
account the cases $\a=\b$ and $\a<\b$.  In the case of logic \iel, rules handling $\K$
have to take into account the relation \E, which is a subset of $\leq$. If $\K A$ is
satisfied by $\a$, to draw a correct deduction we need to know if $\E(\a,\a)$ holds. Thus
calculus \Liel uses the sequent $\eseq{\T}{\G}{\D}$. We say that $\a$ satisfies
$\eseq{\T}{\G}{\D}$, and we write $\a\sat\eseq{\T}{\G}{\D}$ iff $\E(\a,\a)$ and
$\a\sat\seq{\T}{\G}{\D}$ hold. Thus $\eseq{\T}{\G}{\D}$ stores the information that it is
satisfied in a world that \E-reaches itself. if $\a\sat\eseq{\T}{\G}{\D}$ holds, then we
conclude that $\E(\a,\a)$ holds; if $\a\sat \seq{\T}{\G}{\D}$ holds, then we cannot draw
any conclusion about $\E(\a,\a)$.

To define the deduction in our calculi, we need to identify a particular type of sequents
the we call {\em terminal}. We divide terminal sequents in two disjoint categories: {\em
  axioms} and {\em flat}. For sake of concreteness, in the case of calculus \Liel of
Figures~\ref{fig:lsj} and~\ref{fig:rules-for-k}, the axioms are of the kind ~$\ruleIRR$,
$\ruleID$, $\ruleeIRR$ or~$\ruleeID$ and the flat sequents fulfils the following
conditions: $\G\subseteq\PV$, $\D\subseteq \AT$ and $\G\cap\D=\emptyset$, where $\G$ is
the second compartment and $\D$ is the third.

As regards proof construction, we use the rules bottom-up, thus by {\em instantiating the rule
  $\rmR$ with a sequent $\s$} we mean that $\s$ is not a terminal sequent and we use $\s$
to instantiate the conclusion of~$\rmR$. We call {\em the result of the instantiation of
  $\rmR$ with $\s$} the sequents $\s_1,\dots,\s_r$ occurring in the premise of
$\rmR$. We consider the sequents occurring in the premise of the rules enumerated from
left to right, thus $\s_1$ is the leftmost premise of $\rmR$ and $\s_r$ is the rightmost.
We assume that in the instantiation the formulas in evidence in a compartment do not occur
in the set in evidence in the same compartment.



Let $\calC$ be a sequent calculus presented in this paper. Given a sequent $\s$, a {\em
  \mbox{($\calC$-)}tree $\calT$ of sequents for $\s$} fulfils the following properties: (i) the root of
$\calT$ is $\s$; (ii) for every sequent $\s_0$ occurring in $\calT$ as non-leaf node, if
$\s_1,\dots,\s_n$ are the children of $\s_0$, then there exists a rule
$\rmR$ of $\calC$ such that if $\s_0$ instantiates $\rmR$, then $\s_1,\dots,\s_n$ is
the result of the instantiation of $\rmR$ with $\s_0$ and the  sequents
$\s_1,\dots,\s_n$ are enumerated considering from left to right the  sequents in
the premise of $\rmR$.
%
%
$\calT$ is a {\em completed (tree of
  sequents)} if the leaves of $\calT$ are terminal sequents. $\calT$ is a {\em
  $\calC$-proof of $\s$ } if all the leaves are axiom sequents. In this case we say that
{\em $\s$ is provable (in $\calC$)} or {\em $\calC$ proves $\s$}.


\section{A calculus to prove the validity of \iel}
\label{sec:iel}
In this section we discuss the problem of proving the validity in \iel of a given formula
$A\in\calL$ by means of the calculus \Liel provided in Figures~\ref{fig:lsj}
and~\ref{fig:rules-for-k}. To decide the validity of $A$ we look for a (\Liel-)proof of
$\seq{}{}{A}$. If such a proof exists we say that {\em $A$ is provable in \Liel} or {\rm
  Liel proves $A$.}

The rules of \Liel in Figure~\ref{fig:lsj} characterise propositional intuitionistic logic
and are discussed in~\cite{FerFioFio:2013}. Thus we only discuss the rules in
Figure~\ref{fig:rules-for-k}.

By inspection of the rules it is easy to prove that the depth of every tree of sequent
$\calT$ for $\s$ is bounded by the number of connectives occurring in $\s$. As a matter of
fact for every rule $\rmR$, the number of connectives of the sequents occurring in the
premise of $\rmR$ is greater than the number of connectives occurring in the sequent of
the conclusion. From this it follows  that the length of every branch of $\calT$ is
bounded by the number of connectives of $\s$.

The correctness of the rules in Figure~\ref{fig:rules-for-k} is based on the Kripke
semantics of the connective $\K$.
\begin{theorem}[Correctness]
  Let $\uK=\langle S, \rho, \leq,\E, \forza\rangle$ a Kripke model for \iel and $\a\in
  S$. For every rule $\rmR$ of \Liel, if $\a$ satisfies the sequent in the conclusion of
  $\rmR$, then $\a$ satisfies at least a sequent in the premise.
\end{theorem}
{\em Proof. }
We only provide the cases related to \iel.  The correctness of the rules for $\K$ exploits
the fact that from properties~\ref{Im1}-\ref{Im3}, follows that for every $\a\in S$, there
exists $\b\in S$ such that $\a\leq\b$ and $\E(\a,\b)$.
\medskip\\
Rule ($e\K L$): let us suppose that $\a\in S$ satisfies the sequent in the bottom of rule
$e\K L$. Thus, by definition of \E-sequent, $\E(\a,\a)$ holds and by semantical definition
of $\K$, $\a\forza A$ holds. Thus we have proved that $\a$ satisfies the \E-sequent in the
premise of $(e\K L)$.
\medskip\\
Rule ($\K L$): let us suppose that $\a\in S$ satisfies the sequent in the bottom of rule
($\K L$). We recall that by definition of satisfiability of a sequent, $\a$ satisfies the
formulas in the second compartment and for every $\b\in S$, if $\a<\b$, then
$\b\forza \bigwedge\T$.  Moreover, by semantical definition of $\K$, for every $\K B$ in
the second compartment and for every $\b\in S$, if $\E(\a,\b)$ holds, then $\b\forza
B$. Let us consider a final world $\g$ of $\uK$ such that $\a\leq\g$.
Since $\g$ is a final world, by \iel semantics it follows that $\E(\g,\g)$
holds, because on the final worlds relation $\E$ is reflexive. 
We have two cases:
(i) $\a=\g$. Thus $\a$ has no any immediate successor and $\a$ satisfies the leftmost
premise of~($\K L$); (ii) $\a<\g$. Thus $\g$ satisfies the rightmost premise of ($\K L$).
\medskip
\\
Rule ($\K R$): let us suppose that $\a\in S$ satisfies the sequent in the bottom of
rule~($\K R$). This implies that $\a$ forces the formulas in the second compartment and
for every $\b\in S$, if $\a<\b$, then $\b\forza \bigwedge\T$.  Moreover, $\a$ does not
force any formula in the third compartment. Since by hypothesis $\a\nonforza \K B$, there
exists a world $\b\in S$ such that $\a\leq\b$, $\E(\a,\b)$ and $\b\nonforza B$. By the
semantics of connective $\K$, we have that if $\a\forza \K A$, then $\b\forza A$. Now, on
$\b$ we have two cases: (i) $\a=\b$. Then $\a$ satisfies sequent
$\seq{\T}{\G,A_1,\dots, A_n}{B,\D}$ and $\E(\a,\a)$ holds. Thus the leftmost premise is
satisfied; (ii) $\a<\b$. Then $\b$ satisfies sequent $\seq{}{\T, \G,
  A_1,\dots,A_n}{B}$. Thus the rightmost premise is satisfied.  \medskip
\\
Rule ($e\K R$): let us suppose that $\a\in S$ satisfies the sequent in the bottom of
rule~($e\K R$). This implies that $\E(\a,\a)$ holds and $\a$ does not force any formula in
the third compartment. Thus there there exists a world $\b\in S$ such that $\a\leq\b$,
$\E(\a,\b)$ and $\b\nonforza B$. By the property of forcing relation, $\a\nonforza B$
follows. Since $\E(\a,\a)$ holds, also $\a\sat \eseq{\T}{\G}{B,\D}$ holds.  \medskip
\\
Rule ($e\to L$): let us suppose that $\a\in S$ satisfies the sequent in the bottom of the
rule~($e\to L$). Thus $\E(\a,\a)$ holds. By definition of intuitionistic implication at
least one of the following points holds: (i) $\a\forza B$. Thus $\a$ satisfies the
leftmost premise of $e\to L$; (ii) $\a\nonforza A$ and for every $\b\in S$, if $\a<\b$,
then $\b\forza A$. Thus $\b\forza B$. Moreover by the hypothesis and $\a<\b$ follows
$\beta\forza \bigwedge\T$. Hence we have proved that $\a$ satisfies the second premise of
$e\to L$; (iii) there exists $\b\in S$ such that $\a<\b$, $\b\nonforza A$ and for every
$\g\in S$, if $\b<\g$, then $\g\forza A$. Then $\g\forza B$.  Moreover, from $\a<\beta$
and the hypothesis, follows that $\beta\forza\bigwedge\T$.  Thus we have proved that $\b$
satisfies the rightmost premise of~($e\to L$).  \qed

We remark that in the proof of the correctness of rule~($\K R$), we do not claim that $\b$
is a final world. This explains the occurrence $\emptyset$ in the second premise of~($\K R$).
The presence of $\bot$ in the first compartment of the premises of~($\K L$) expresses the
fact that the sequents, if satisfiable, must be satisfied by Kripke models containing
exactly one world. The effect in rule application is that when this kind of sequent
instantiates the rules ($\to R$), ($e\to R$), ($\to L$), ($e\to L$), ($\K R$) or~($\K L$), the
resulting $\s_r$ (the rightmost premise) is the axiom $\ruleIRR$ or $\ruleeIRR$.

\begin{example}\rm
  Calculus \Liel proves~\ref{Ax2} $\K(A\to B)\to (\K A\to \K B)$.  We remark that when
  $\T=\emptyset$ and $\D=\bot$ the two premises of rule~($\to R$) coincide, thus we only show one
  branch.
  \newcommand{\seqa}{\seq{}{}{\K(A\to B)\to (\K A\to \K B)}}
  \newcommand{\seqb}{\seq{}{\K(A\to B)}{ \K A\to \K B}} \newcommand{\seqc}{\seq{}{\K(A\to
      B),\K A}{\K B}}
  \newcommand{\seqd}{\seq{}{A\to B, A}{B}}
  \newcommand{\eseqd}{\eseq{}{A\to B, A}{B}}
  \newcommand{\axa}{\seq{}{B, A}{B}~\ruleID}
  \newcommand{\axb}{\seq{B}{ A}{A, B}~\ruleID}
  \newcommand{\axc}{\seq{B}{A}{A}~\ruleID}
  \newcommand{\eaxa}{\eseq{}{B, A}{B}~\ruleID}
  \newcommand{\eaxb}{\eseq{B}{ A}{A, B}~\ruleID}
  \newcommand{\deda}{\regola{\seqd}{\axa & \axb & \axc}{\to L}}
  \newcommand{\ededa}{\regola{\eseqd}{\eaxa & \eaxb & \axc}{e\to L}}
  \newcommand{\dedb}{\regola{\seqc}{\lbrack\textrm{Ded}_1\rbrack & \lbrack\textrm{Ded}_2\rbrack}{\K R}}
  \newcommand{\dedc}{\regola{\seqb}{\dedb}{\to R}}
  \newcommand{\dedd}{\regola{\seqa}{\dedc}{\to R}}
  \newcommand{\dedxx}{\apply \seqb \seqa {\to R}}
  \[
    \begin{array}[t]{c}
      \begin{array}[c]{cc}
        \lbrack\textrm{Ded}_1\rbrack
        &
          \begin{array}{c}
            \ededa
          \end{array}
          \medskip\\
        \lbrack\textrm{Ded}_2\rbrack
        &
          \begin{array}{c}
            \deda
          \end{array}
      \end{array}
      \bigskip\\
      \dedd
    \end{array}
    %
  \]
\end{example}

\begin{example}
  \rm Calculus \Liel proves~\ref{Ax4} $\K A\to \neg\neg A$.
  Rules~($\K L$) and~($\to R$) are instantiated with $\D=\{\bot\}$ and $\T=\emptyset$. Thus
  the premises are equal and to save space we only show one branch.
  \newcommand{\axa}{\eseq{}{A,\bot}{\bot}~\ruleIRR}
  \newcommand{\axb}{\eseq{\bot}{A}{A,\bot}~\ruleID}
  \newcommand{\axc}{\eseq{\bot}{A}{A}~\ruleID}
  \newcommand{\seqa}{\eseq{}{A,\neg A}{}}
  \newcommand{\deda}{\regola{\seqa}{\axa & \axb & \axc}{e\to L}}
  \newcommand{\seqb}{\seq{}{\K  A,\neg A}{}}
  \newcommand{\dedb}{\regola{\seqb}{\deda}{\K L}}
  \newcommand{\seqc}{\seq{}{\K A}{\neg\neg A}}
  \newcommand{\dedc}{\regola{\seqc}{\dedb}{\to  R}}
  \newcommand{\seqd}{\seq{}{}{\K A\to\neg\neg A}}
  \newcommand{\dedd}{\regola{\seqd}{\dedc}{\to R}}
  \[
    \dedd
  \]
\end{example}

\begin{example}\rm
Calculus \Liel proves the double negation of the classical reflection: $\neg\neg(\K A\to
A)$. The instantiation of ($\to R$) gives two equal premises, thus we do not show the
branch. Moreover, since $\bot$ in the third compartment is irrelevant, we disregard it in
the second premise of ($\to L$) application. 
\newcommand{\eaxa}{\eseq{\bot}{A}{A}~\ruleID}
\newcommand{\axb}{\seq{}{\bot,A}{\bot}~\ruleIRR}
\newcommand{\eaxb}{\eseq{}{\bot,A}{\bot}~\ruleIRR}
\newcommand{\axd}{\seq{}{\bot,\K A}{A}~\ruleIRR}
\newcommand{\seqa}{\seq{\bot}{\K A}{A}}
\newcommand{\deda}{\regola{\seqa}{\eaxa & \eaxb}{\K L}}
\newcommand{\seqb}{\seq{\bot}{}{\K A\to A}}
\newcommand{\dedb}{\regola{\seqb}{\deda & \axd}{\to R}}
\newcommand{\seqc}{\seq{}{\neg(\K A\to A)}{\bot}}
\newcommand{\axc}{\seq{}{\false}{\false}~\ruleIRR}
\newcommand{\dedc}{\regola{\seqc}{\axc & [\textrm{Ded}] & [\textrm{Ded}]}{\to L}}
\newcommand{\seqd}{\seq{}{}{\neg\neg(\K A\to A)}}
\[
  \begin{array}[c]{cc}
    [\textrm{Ded}] &
          \begin{array}{c}\footnotesize
            \dedb\medskip
          \end{array}\medskip\\
        &\regola{\seqd}{\dedc}{\to R}
  \end{array}
\]
\end{example}

\begin{example}\rm
  \Liel does not prove the classical reflection $\K A\to A$. As a matter of fact, the
  following is the unique sequent tree for $\K A\to A$ that we can build with the rules of
  \Liel.
  \newcommand{\seqa}{\eseq{}{A}{}}
  \newcommand{\axa}{\eseq{}{A}{A}\ruleID}
  \newcommand{\seqb}{\seq{}{\K A}{A}}
  \newcommand{\seqc}{\seq{}{}{\K A\to A}}
  \newcommand{\deda}{\regola{\seqb}{\axa & \seqa}{\K L}}
  \[
    \regola{\seqc}{\deda}{\to R}
  \]
  In the following we provide a procedure that failing to return a proof, returns a Kripke
  model that satisfies the $\seqc$, that is a model that does not force $\K A\to A$. 
\end{example}

\section{Completeness}
\label{sec:liel:completeness}
In the following we design a procedure that, given a sequent, uses \Liel to build a proof,
if any, otherwise returns a model whose root satisfies the sequent.
A feature of \Liel is
that if all the branches generated by rule instantiation are developed in systematic way,
then if a model exists, the failure in proof search allows to get a model of minimum depth.
This result is attained by the way \Liel handles the connectives $\to$ and $\K$.
Our procedure is designed to return models of minimum depth.
\medskip

%

\noindent
{\sc Procedure $\proc(\s)$} 
\begin{enumerate}[label=(\arabic*),ref=\arabic*]
  
\item \label{assiomi}\label{proc:assiomi}%
  if $\s$ can instantiate an axiom, then return the proof \[\s~(\rmR)\] where
  $\rmR\in\{\mathbf{Irr}, \mathbf{Id},\mathbf{eIrr},\mathbf{eId}\}$;

\item\label{no-rule-to-apply}\label{proc:no-rule-to-apply}%
  if $\s$ is a flat sequent $\seq{\T}{\G}{\D}$ or $\eseq{\T}{\G}{\D}$, then return the
  structure
  \[
    \uK=\langle  S, \rho, \leq, \E, \forza\rangle
  \]
  where $S$ = $\{\rho\}$, $\leq$ = $\{(\rho,\rho)\}$, $\E$ = $(\rho,\rho)$ and $\forza$ =
  $\{\rho\}\times \G$;
 
\item \label{proc:alpha-step} if $\s$ can instantiate $\rmR\in\{\land L, \lor R,e\land L,
  e\lor R, e\K L, e\K R \}$, then let
  $\s_1$ be the result of an instantiation of $\rmR$ with $\s$. Let $U=\proc(\s_1)$. If
  $U$ is a structure, then return $U$, otherwise return $\apply{U}{\s}{(\rmR)}$;
  
\item \label{proc:beta-step} if $\s$ can instantiate $\rmR\in\{\lor L, \land R, e\lor L,
  e\land R \}$, then let
  $\s_1$ and $\s_2$ be the result of an instantiation of $\rmR$ with $\s$.
  For $i=1,2$, let $U_i=\proc(\s_i)$;
  \begin{enumerate}
  \item If $U_1$ and $U_2$ are proofs, then return $\regola{\s}{U_1 & U_2}{\rmR}$;
  \item \label{proc:beta-step:exactly-one} if exactly one between $U_1$ and $U_2$ is a
    structure, then return the structure $U_i$, where $i\in\{1,2\}$;
  \item
    \label{proc:beta-step:both}
    If $\depth(U_1)<\depth(U_2)$, then return $U_1$, else return $U_2$;
  \end{enumerate}
\item   \label{proc:non-invertile-step}
  if $\s$ can instantiate one of the rules in $\{(\to R)$, ($\to L$), ($\K R$), ($e\to R$), ($e\to L)\}$ then:
  \begin{enumerate}[label=(\alph*),ref=\theenumi.\alph*]
  \item let $\noninv=\emptyset$ and $\inv=\emptyset$;
  \item
    \label{proc:non-invertile-step:loop}
    for every $\rmR\in \{\to R, \K R, e\to R\}$, for every possible instantiation of
    $\rmR$ with $\s$:
    \begin{enumerate}[label=(\roman*),ref=\theenumii.\roman*]
    \item \label{proc:non-invertile-step:loop:calls}
      let $\s_1$ and $\s_2$ the result of the instantiation of $\rmR$ with $\s$, let
      $U_1=\proc(\s_1)$ and $U_2=\proc(\s_2)$;
    \item \label{proc:proof-step} if $U_1$ and $U_2$ are proofs, then return
      $\regola{\s}{U_1 & U_2}{\rmR}$;
    \item \label{proc:non-invertile-step:loop:first-update-U1}
      if $U_1$ is a structure, then let $\inv=\inv\cup\{U_1\}$;
    \item
      \label{proc:non-invertile-step:loop:first-update-U2}
      if $U_2$ is a structure, then let $\noninv=\noninv\cup\{U_2\}$
    \end{enumerate}
  \item for every $\rmR\in \{\to L, e\to L\}$, for every possible instantiation of
    $\rmR$ with $\s$:    
    \begin{enumerate}[label=(\roman*),ref=(\roman*)]
    \item let $\s_1$, $\s_2$ and $\s_3$ the result of the instantiation and let 
      $U_1=\proc(\s_1)$, $U_2=\proc(\s_2)$ and $U_3=\proc(\s_3)$;
    \item if $U_1$, $U_2$ and $U_3$ are proofs, then return $\regola{s}{U_1 & U_2 &U_3}{\to
        L}$;
    \item if $U_1$ is a structure, then let $\inv=\inv\cup\{U_1\}$;
    \item if $U_2$ is a structure, then let $\inv=\inv\cup\{U_2\}$;
    \item if $U_3$ is a structure, then let $\noninv=\noninv\cup\{U_3\}$;
      
    \end{enumerate}

  \item \label{proc:glue-step} Let $\uM_1,\dots,\uM_n$ be an enumeration of the structures
    in $\noninv$, where for $i=1,\dots,n$,
    $\uM_i=\langle S_i,\rho_i,\E_i,\leq_i,\forza_i\rangle$. We suppose that for
    $i,j=1,\dots,n$, if $i\neq j$, then $S_i\cap S_j=\emptyset$. Let
    $\uM=\langle S,\rho,\E,\leq,\forza\rangle$ be defined as follows:
    \begin{itemize}
    \item[-] $\displaystyle S=\{\rho\}\cup\bigcup_{i=1}^n S_i$;
      
    \item[-]
      $\displaystyle \leq=\Big\{(\rho,\alpha)\in \{\rho\}\times S\Big\}\cup\bigcup_{i=1}^n
      \leq_i$;
      
    \item[-] $\setlength{\arraycolsep}{2pt}\begin{array}[t]{lll}
        \textstyle\E' = \bigcup_{i=1}^n \E_i\cup \{(\rho,\a)\in \{\rho\}\times S
        & | &
              \textrm{there exits $\b\in S$ such that }\\
        & & (\b,\a)\in \bigcup_{i=1}^n \E_i \}
      \end{array}
      $ 
    \item[] if $\s=\eseq{\T}{\G}{\D}$, then $\E=\E'\cup\{(\rho,\rho)\}$ else $\E=\E'$;
      
    \item[-] $\displaystyle \forza\,\,=\,\Big\{(\rho, p)\in \{\rho\}\times \PV\ |\ p\in
      \G\Big\}\cup\bigcup_{i=1}^n \forza_i$
    \end{itemize}
  \item
    \label{proc:glue-model}
    if $\inv=\emptyset$, then return $\uM$;
  \item let $\mindepth = \min \{\depth(U) | U\in\inv \}$ and\\ let
    $U\in\{U'\in\inv\, |\, \depth(U')=\mindepth\}$;
  \item
    \label{proc:inv-model}
    if the cardinality of $\noninv$ is less than the number of possible
    instantiations of the rules~($\to L$), ($\to R$), ($\K R$), ($e\to R$) and ($e\to L$) with
    $\s$, then return $U$;
  \item \label{proc:comparison}
    if $\depth(\uM)< \mindepth$, then return $\uM$, else return $U$; 
      
  \end{enumerate}
  
\item
  \label{proc:rule-K-L}
  let $\s_1$ and $\s_2$ be the result of an
  instantiation of rule ($\K L$) with $\s$. Let $U_i=\proc(\s_i)$, for $i=1,2$;
  \begin{enumerate}
  \item If $U_1$ and $U_2$ are proofs, then return $\regola{\s}{U_1& U_2}{\K L}$;
  \item \label{proc:KL-return-invertible}%
    if $U_1$ is a structure and $U_2$ a proof, then return $U_1$;
  \item
    \label{proc:rule-K-L:build-structure}
    Let $U_2$ be the structure $\langle S',\rho', \E', \leq',\forza'
    \rangle$. Let $\uM$ be the structure $\langle S,\rho,\E,\leq,\forza\rangle$ defined as
    follows:
    \begin{itemize}
    \item[-] $\displaystyle S=\{\rho\}\cup S'$;
      
    \item[-] $\displaystyle \leq=\Big\{(\rho,\alpha)\in \{\rho\}\times S\Big\}\cup\leq'$;
      
    \item[-] 
      $
        \displaystyle\E=\Big\{ (\rho,\alpha)\in \{\rho\}\times S\, |\, (\rho',\alpha)\in\E'\Big\}\,
        \cup\E';
      $
      
    \item[-] $\displaystyle \forza\,\,=\,\Big\{(\rho, p)\in \{\rho\}\times \PV\ |\ p\in
      \G\Big\}\,\cup\forza'$.
    \end{itemize}
  \item \label{proc:KL-return-glued} \label{proc:rule-K-L:return-structure}%
    if $U_1$ is a proof, then return $\uM$;
  \item \label{proc:KL-return-comparison}%
    if $\depth(U_1)<\depth(\uM)$, then return $U_1$,
    else $\uM$.
  \end{enumerate}
\end{enumerate}
{\sc End Procedure \proc}\medskip\\
We remark that if Step~\ref{proc:rule-K-L} is reached, then we know that $\s$ instantiates
rule $\K L$.

Note  that the satisfiability of the rightmost premise of $\K L$ does not imply the
satisfiability of the sequent in the conclusion. In other words, the unprovability of the
rightmost premise does not imply the unprovability of the sequent in conclusion. Thus
backtracking can be required.  In the completeness theorem we prove that, when the third
compartment contains atoms only and the second compartment contains propositional
variables and $\K$-formulas only, then from the satisfiability of the rightmost premise we
deduce the satisfiability of the conclusion. In other terms, under the stated conditions,
if rule $\K L$ is instantiated, then the backtracking can be avoided.

We start by proving that if \proc returns a structure $\uK$, then $\uK$ fulfils the
definition of Kripke model for \iel. This is the first point to prove the completeness theorem.

\begin{theorem}[Kripke Model]
  \label{th:kripke-model}
  \rm Let $\s$ be a sequent . If $\proc(\s)$ returns a structure $\uK$, then $\uK$ is a
  Kripke model for \Liel.
\end{theorem}
{\em Proof.}
We proceed by induction on the depth of the recursive calls.
\medskip\\
{\em Base:}
$\proc(\s)$ returns $\uK$ without performing any recursive call. Then $\uK$ is the result
of the construction in Step~\ref{no-rule-to-apply} and $\s$ is flat.  It is immediate to
check that the definition of $\uK$ at Step~\ref{no-rule-to-apply} fulfils the definition
of Kripke model for \iel. \medskip\\
{\em Induction:}
\proc returns $\uK$ by performing recursive calls that by induction hypothesis return a
Kripke model.  We proceed by considering every step of \proc where a structure is
returned.
\medskip\\
{\bf $\uK$ is the result of performing Step~(\ref{proc:alpha-step}),
Step~\ref{proc:beta-step:exactly-one} or Step~\ref{proc:beta-step:both}.} Immediate, since the
elements returned by \proc are built by the recursive calls and hence, by induction
hypothesis, they are Kripke models.\medskip\\   
%
%
{\bf $\uK$ is the result of performing Step~(\ref{proc:glue-model}.}) We have to show that
the returned structure $\calM$ is a Kripke model. Here we have two main cases: (i) $\s$ is
an \E-sequent. The rules that can be instantiated with $\s$ are ($e\to L$) and ($e\to R$);
(ii) $\s$ is not an \E-sequent. The rules that can be instantiated with $\s$ are ($\to L$)
($\to R$) and $(\K R)$. In the following we handle the two cases in a row, without any
further distinction.  The sequents $\s_1$ and $\s_2$ resulting from every possible
instantiation of the rules $\{(\to L), (e\to L)\}$ with $\s$ and the sequent $\s_1$ resulting
from every possible instantiation of the rules $\{(\to R), (e\to R), (\K R) \}$ with $\s$ are
provable, that is there is no Kripke model satisfying them. By construction of \proc,
follows that if Step~\ref{proc:glue-model} is performed, then $\s_r$ resulting from every
possible instantiation of the rules $\{\to L, e\to L,\to R, e\to R, \K R \}$ with $\s$ has
a Kripke model (otherwise a proof would have been returned). Procedure \proc collects all
these models in $\noninv$.
Since for every element $\a$ of $\calM$, $\rho\leq\a$ holds, for every $U\in \noninv$,
the root of $U$ forces $\bigwedge\T$ and $\bigwedge\G$ and by definition of forcing in
$\rho$, follows that the returned structure $\uM$ obeys to the definition of Kripke model
for propositional intuitionistic logic.
As regards the definition of \E, for every
$U\in\noninv$, if $\E'(\rho',\a)$ holds, where $\rho'$ is the root of $U$, $\a$ a world of
$U$ and $\E'$ the $\E$-relation of $U$, then, $\E(\rho, \a)$ holds by construction of
$\uM$. Thus $\uM$ obeys to the definition of Kripke model for \iel and we have proved the
statement of the theorem.
\medskip\\
{\bf $\uK$ is the result of performing Step~(\ref{proc:inv-model}),
Step~(\ref{proc:comparison}) or Step~\ref{proc:KL-return-invertible}.} Immediate by using
induction hypothesis.
\medskip\\
{\bf $\uK$ is the result of performing Step~\ref{proc:KL-return-glued}.} Analogous to the proof
for Step~\ref{proc:glue-model} applied to only one model (that is $n=1$).

{\bf $\uK$ is the result of performing~\ref{proc:KL-return-comparison}.}
Immediate since $U_1$ and $\calM$ are structures that \proc
returns in Step~\ref{proc:KL-return-invertible} or Step~\ref{proc:KL-return-glued} and we
have already proved that they satisfy the statement of the theorem.
\qed


Now that we have proved that given a sequent $\s$, if $\proc(\s)$ returns a structure
$\uK$, then $\uK$ is a Kripke model for \iel, in the following we prove that the structure $\uK$
satisfies $\s$. This is the main point to prove the completeness. 
\begin{proposition}[Satisfiability]
  \label{prop:satisfiability}
  \rm Let $\s$ be a sequent . If $\proc(\s)$ returns a structure $\uK$, then the root of
  $\uK$ satisfies $\s$.
\end{proposition}
{\em Proof.} By the rules of the calculus, the sequent $\s$ is of the kind
\mbox{$\seq{\T}{\G}{\D}$} or \mbox{$\eseq \T \G \D$}. We recall from
Section~\ref{sec:definitions} that if $\eseq{\T}{\G}{\D}$ is satisfied, then
$\seq{\T}{\G}{\D}$ is satisfied.  We proceed by induction on the depth of the recursive
calls.
\medskip\\
{\em Base:}
$\proc(\s)$ returns $\uK$ without performing any recursive call. Then $\uK$ is the result
of the construction in Step~\ref{no-rule-to-apply} and $\s$ is flat.  Since $\s$ is flat,
from $\G\subseteq \PV$ and by definition of $\forza$, immediately follows that
$\rho\forza\bigwedge\G$ holds; from $\D\subseteq \AT$ and $\G\cap \D=\emptyset$ follows that for
every $B\in\D$, $\rho\nonforza B$. By the fact that $\uK$ only contains a world, $\rho$
satisfies the first compartment of $\s$ and $\E(\rho,\rho)$ holds. Thus we have proved
that
$\rho\sat \s$.\medskip\\
{\em Induction:} \proc returns $\uK$ by performing recursive calls that by induction
hypothesis return Kripke models fulfilling the statement of the proposition.  We proceed by
considering every step of \proc where a Kripke model is returned.
\medskip\\
%
%
{\bf $\uK$ is the result of performing Step~(\ref{proc:alpha-step}).} By induction hypothesis  $\uK$
satisfies $\s_1$, where $\s_1$ is the result of the instantiation of
$\rmR\in\{(\land L), (\lor R), (e\land L), (e\lor R), (e\K L), (e\K R) \}$ with $\s$. We proceed by
cases:
if $\rmR$ is~$(e\K R)$, then $\s_1$ is of the kind \mbox{$\eseq{\T}{\G}{B,\D}$}. Since by
induction hypothesis $\uK$ satisfies $\s_1$, we have that $\rho\nonforza B$ and
$\E(\rho,\rho)$ hold. By the semantics of $\K$ follows that $\rho\nonforza \K B$ and thus
$\rho$ satisfies $\s$;
if $\rmR$ is ($e\K L$), then $\s_1$ is of the kind $\eseq{\T}{\G, A}{\D}$. Since
$\rho\forza A$, for the persistence property, $A$ holds in every world of $\uK$ an we
immediately get $\rho\forza \K A$ and $\rho$ satisfies $\s$.
The other cases are similar.\medskip\\
%
%
{\bf $\uK$ is the result of performing Step~\ref{proc:beta-step:exactly-one}.} By
induction hypothesis, $\uK$ satisfies $\s_i$. To prove that $\uK$ satisfies $\s$ we
proceed by cases according to the possible values of $\rmR$. If $\rmR$ is ($e\lor L$), then,
since by induction hypothesis $\uK$ satisfies $\s_i$, it follows that $\E(\rho,\rho)$ and
$\rho\forza A\lor B$ hold. Thus we
have proved that $\uK$ satisfies $\s$. The other cases of $\rmR$ are similar.\medskip\\
%
{\bf $\uK$ is the result of performing Step~\ref{proc:beta-step:both}.} By induction
hypothesis $U_1$ and $U_2$ fulfil the statement of the proposition on respectively
$\s_1$ and $\s_2$. At this point we get the statement of the proposition on $\s$ by applying to the
returned model the proof provided for Step~\ref{proc:beta-step:exactly-one}.
\medskip\\
%
%
{\bf $\uK$ is the result of performing Step~(\ref{proc:glue-model}).} At this stage of the
procedure, the formulas in $\G\cup\D$ are atoms, implications and $\K$-formulas. We remark
that if $\s=\eseq{\T}{\G}{\D}$, then no $\K$-formula belongs to $\G$.  We prove that
$\rho$ satisfies every compartment of $\s$:
\begin{itemize}[label=-]
\item
  for every $A\in\D\cap\at$, $\rho\nonforza A$ follows from 
  the definition of $\forza$, $\D\subseteq\at$ and
  $\G\cap\D=\emptyset$;
\item
  every root of every model in $\noninv$ forces the formulas in the first compartment of $\s$;
\item
  for every $A\to B\in \G$, let $\s_3$ be the result of instantiating $\s$ with ($\to L$) or
  ($e\to L$) according to the type of $\s$. The recursive call $\proc(\s_3)$ returns a model
  $U$ that is collected in $\noninv$. Let $\rho'$ be the root of $U$.  By
  induction hypothesis $\rho'\sat\s_3$, hence $\rho'\nonforza A$ and for every world
  $\a$ of $U$ different from $\rho'$, $\a\forza B$. By the meaning of implication, we
  have that $\rho'\forza A\to B$. We remark that all the $\K$-formulas in $\G$ are in the
  second compartment of $\s_3$, thus, by induction hypothesis, for every $\K B\in\G$,
  $\rho'\forza \K B$.  By construction of $\calM$ and
  Theorem~\ref{th:kripke-model}, $\rho\nonforza A$. For every model
  $U'\in\noninv$ different from $U$, we have that $U'=\proc(\s')$ and $A\to B$ is
  in the second compartment of $\s'$. By induction hypothesis $U'$ satisfies $\s'$ and
  thus the root of $U'$ forces $A\to B$. Summarising, we have proved that for every
  $A\to B\in\G$, $\rho\forza A\to B$;
\item
  for every $A\to B\in\D$, the recursive call $\proc(\s_2)$ returns a model $U$
  collected in $\noninv$, with $\s_2=\seq{}{\G,\T,A}{B}$.  Let $\rho'$ be the root of
  $U$.  By induction hypothesis, $\rho'\sat\s_2$, hence $\rho'\forza A$ and
  $\rho'\nonforza B$. Thus $\rho'\nonforza A\to B$.  By construction of $\calM$ and
  Theorem~\ref{th:kripke-model} we have that $\rho\nonforza A\to B$. Thus we
  have proved that for every $A\to B\in\D$, $\rho\nonforza A\to B$;

\item %
  for every $\K B\in \D$, the recursive call $\proc(\s_2)$ returns a model $U$
  collected in $\noninv$, with $\s_2=\seq{}{\G,\T, A_1,\dots,A_n}{B}$. By induction
  hypothesis, $\rho'\forza A_1,\dots,\rho'\forza A_n$, thus
  $\rho'\forza \K A_1,\dots,\rho'\forza \K A_n$, and $\rho'\nonforza B$, with $\rho'$ root
  of $U$.  By construction of $\calM$, $\E(\rho,\rho')$ holds. We get
  $\rho\nonforza \K B$.  Thus we have proved that for every $\K B\in\D$,
  $\rho\nonforza \K B$.

\end{itemize}
In the points above we have proved that for every $\K A\in\G$ and for every model
$U\in\noninv$, the root of $U$ forces $\K A$. Moreover by last point above, if $E(\rho,\rho')$
holds, then $\rho'\forza A$. Thus we have proved that for every $\K A\in\G$,
$\rho\forza \K A$.  Summarising the points above prove that $\rho$ satisfies
$\s=\seq{\T}{\G}{\D}$.  Finally, by construction of $\calM$, if
$\s=\eseq{\T}{\G}{\D}$, $\E(\rho,\rho)$ holds and thus we have proved that
$\rho\sat \s$ for any kind of $\s$.
\medskip\\
{\bf $\uK$ is the result of performing
  Step~(\ref{proc:inv-model}).} The point is proved by cases. We provide two of them:
\begin{itemize}
\item $U=\proc(\s_1)$, where $\s_1$ is the result of instantiating ($\K R$) with
  $\s$. Since $U=\uK$, by induction hypothesis the root $\rho$ of $\uK$ satisfies
  $\s_1$. This implies that $\rho\nonforza B$, $\rho\forza A_1,\dots,\rho\forza A_n$ and
  $\E(\rho,\rho)$ hold. We immediately get that  $\rho\forza
  \K A_1,\dots,\rho\forza \K A_n$ and $\rho\nonforza \K B$ hold. Thus we have proved that
  $\rho$ satisfies $\s$; 
\item $U=\proc(\s_2)$, where $\s_2$ is one of the sequent in the result of the
  instantiation of ($e\to L$) with $\s$. Since $U=\uK$, by induction hypothesis the
  following hold: $\rho\forza\bigwedge\T$; $\E(\rho,\rho)$; $\rho\nonforza A$; for every
  $\a\in S$, if $\rho<\a$, then $\a\forza B$. By the meaning of implication we have proved
  that $\rho\forza A\to B$ and hence $\rho\sat\s$.
\end{itemize}
The other cases are similar.
\medskip\\
%
%
{\bf $\uK$ is the result of performing Step~(\ref{proc:comparison}).} If $\proc(\s)$ returns
$\uM$, then the proof goes as in the case of Step~(\ref{proc:glue-model}), otherwise the
proof goes as in the case of Step~(\ref{proc:inv-model}).
\medskip\\
{\bf $\uK$ is the result of performing Step~\ref{proc:KL-return-invertible}.}
By induction hypothesis, $\uK$ satisfies
$\s_1$. Therefore, $\rho\forza A$ and $\E(\rho,\rho)$ and by the fact that the
first compartment of $\s_1$ contains $\bot$, follows that $S=\{\rho\}$. Thus  $\rho\forza \K
A$ and $\rho$ satisfies the first compartment of $\s$. Thus we have proved that $\rho$
satisfies $\s$.\medskip\\
{\bf $\uK$ is the result of performing Step~\ref{proc:KL-return-glued}.}  By induction
hypothesis $U_2$ is a Kripke model that satisfies $\s_2$. This implies that $S'=\{\rho'\}$ and
$\E'(\rho'\rho')$ hold. $\E'(\rho'\rho')$ implies that for every $\K B\in \G$,
$\rho'\forza B$.  By construction of $\calM$, $\E(\rho,\rho')$ and $\rho\leq\rho'$
hold. Moreover the following holds: for every $\K B\in \G$, $\rho'\forza B$. Thus by the
semantical meaning of $\K$ we have: for every $\K A$ occurring in the second compartment
of $\s$, $\rho\forza \K A$. We recall that if we reach Step~\ref{proc:KL-return-glued} the
second compartment of $\s$ is included in $\PV\cup\{\K A\in\calL\}$ and $\D\subseteq\AT$.
By definition of $\forza$, $\rho$ forces all the propositional variables in the second
compartment of $\s$ and $\rho$ does not force any
element in $\D$. Summarising we have proved that $\rho$ satisfies $\s$.\medskip\\
{\bf $\uK$ is the result of performing~\ref{proc:KL-return-comparison}.}
Immediate since $U_1$ and $\calM$ are structures that \proc
returns in Step~\ref{proc:KL-return-invertible} or Step~\ref{proc:KL-return-glued} and we
have already proved that they satisfy the statement of the proposition.
\qed

From the proposition above we have the completeness theorem:
\begin{theorem}[Completeness]
  \label{th:completeness}
  Let $A\in\calL$. If $\proc(\seq{}{}{A})=\uK$, then $A$ is not valid in \iel. 
\end{theorem}
{\em Proof:} By Proposition~\ref{prop:satisfiability}, the root $\rho$ of $\uK$ satisfies
$\seq{}{}{A}$. By definition of satisfiability, we have $\rho\nonforza A$ and hence $A$ is
not valid in \iel. \qed


%

In the following we show that \proc returns a Kripke model of minimum depth.
\begin{theorem}[Minimality]
\label{th:minimality}
  \rm Let $\s$ be a sequent . If $\proc(\s)$ returns a structure $\uK$, then
  \[
    \depth(\uK)=\min\{\depth(U)\ |\ U\textrm{ is a Kripke model with root $\rho$ and }
    \rho\sat \s\}
    .\]
\end{theorem} {\em Proof.}  We proceed by induction on the depth of the recursive calls.
\medskip\\
{\em Base:}
$\proc(\s)$ returns $\uK$ without performing any recursive call. Then $\uK$ is the result
of the construction in Step~\ref{no-rule-to-apply} and $\s$ is flat.  The statement of the
theorem immediately follows by the
fact that $\uK$ has one world.  \medskip
\\
{\em Induction:} \proc returns $\uK$ by performing recursive calls that by induction
hypothesis return Kripke models fulfilling the statement of the theorem. 
\medskip\\
%
%
{\bf $\uK$ is the result of performing Step~(\ref{proc:alpha-step}).} Let us assume that
there exists a model $\uK'=\langle S',\rho',\leq',\E',\forza'\rangle$ that satisfies $\s$
and $\depth(\uK')<\depth(\uK)$. We use the induction hypothesis that $\uK$ satisfies the
statement of the theorem on $\s_1$ to
get a contradiction. We proceed by cases according to $\rmR$:
if $\rmR$ is ($e\K R$), then $\s$ is of the kind $\eseq{\T}{\G}{\K B,\D}$ and thus
$\rho'\nonforza' \K B$. This implies that there exists $\a\in S'$ such that
$\rho'\leq' \a$, $\E'(\rho',\a)$ and $\a\nonforza' B$. Thus $\rho'\nonforza' B$ and $\uK'$
satisfies $\s_1$. $\uK'$ is a model such that $\depth(\uK')<\depth(\uK)$ and this
contradicts the induction hypothesis that there is no model satisfying $\s_1$ with depth less
than the depth of $\uK$;
if $\rmR$ is ($e\K L$), then $\s$ is of the kind $\eseq{\T}{\G, \K A}{\D}$ and thus
$\rho'\forza' \K A$. Since $\E'(\rho',\rho')$ holds, then $\rho'\forza' A$ holds and thus
$\uK'$ satisfies $\s_1$. This is absurd because by induction hypothesis there is no model
satisfying $\s_1$ with depth less than the depth of $\uK$.
The other cases are similar.\medskip\\
%
%
{\bf $\uK$ is the result of performing Step~\ref{proc:beta-step:exactly-one}.} Let us
assume that there exists $\uK'$ such that $\depth(\uK')< \depth(\uK)$ and $\uK'$ satisfies
$\s$. We go by cases on $\rmR$: if $\rmR$ is ($e\lor L$), then $\uK'\forza' A\lor B$. Moreover
$\E'(\rho',\rho')$ holds. Thus $\uK'$ satisfies $\s_i$ (since $\s_{3-i}$ has a
proof). This contradicts the induction hypothesis on $\uK$ for $\s_i$. The other cases of
$\rmR$ are similar.
\medskip\\
{\bf $\uK$ is the result of performing Step~\ref{proc:beta-step:both}.} By induction
hypothesis $U_1$ and $U_2$ fulfil the statement of the theorem respectively for $\s=\s_1$ and
$\s=\s_2$.  Procedure \proc returns the model of minimum depth between $U_1$ and $U_2$. At
this point we get the statement of the theorem  by applying to the returned model
the proof provided for Step~\ref{proc:beta-step:exactly-one}.
\medskip\\
%
%
{\bf $\uK$ is the result of performing Step~(\ref{proc:glue-model}).} Let
$\md=\max\{\depth(U) | U\in\noninv\}$.  By construction, the depth of $\uK$ is $\md +1$.  Let us suppose
that there exists a model $\uK'$ of depth less than $\md + 1$. By the hypothesis that
$\inv=\emptyset$, it follows:
\begin{enumerate}[label=(\alph*),ref=\alph*]
\item\label{min:a}%
  for every possible instantiation of the rules ($\to L$) and ($e\to L$) with $\s$, the
  sequents $\s_1$ and $\s_2$ are not satisfiable by any Kripke model for \iel. The
  correctness of the rules implies that if $\s$ is satisfied by a world of a model and
  every world of every model does not satisfy $\s_1$ and $\s_2$, then $\s_3$ is
  satisfiable. Sequent $\s_3$ cannot be satisfied by the root of $\uK'$, otherwise,
  together with the satisfiability of $\s$, we would get that the root of $\uK'$ satisfies
  $\s_2$, against the assumption. Thus $\s_3$ is satisfied by a world $\a\in S'$ such that
  $\rho'<'\a$. We conclude that for every possible instantiation of the rules ($\to L$) and
  ($e\to L$) with $\s$, the resulting sequent $\s_3$ is satisfiable by a \iel model whose
  depth is lower than $\md$.
\item\label{min:b}%
  for every possible instantiation of $\s$ with ($\to R), (e\to R$) and ($\K
  R$), the sequent
  $\s_1$ is not satisfiable by any Kripke model for \iel. The correctness of the rules
  implies that $\s_2$ is satisfiable. Sequent
  $\s_2$ cannot be satisfied by the root of
  $\uK'$, otherwise, together with the satisfiability of
  $\s$, we would get that the root of $\uK'$ satisfies
  $\s_1$, against the assumption. Thus $\s_2$ is satisfied by a world $\a\in
  S'$ such that
  $\rho'<'\a$. We conclude that for every possible instantiation of the rules ($\to
  R$), ($e\to R$) and ($\K R$) with $\s$, the resulting sequent
  $\s_2$ is satisfiable by a \iel model whose depth is lower than $\md$.
\end{enumerate}
From Points~(\ref{min:a}) and~(\ref{min:b}) we get that every sequent at hand is
satisfiable with a Kripke model of depth lower than $\md$, but in $\noninv$ there exists a
model whose depth is equal to $\md$. Such a model is returned by a recursive call of
\proc. By induction hypothesis, \proc  returns Kripke models of optimal depth. Thus the model
$\uK'$ does not exist.
\medskip\\
{\bf $\uK$ is the result of performing Step~(\ref{proc:inv-model}).}
Let us suppose there exists a model $\uK'$ such
that $\rho'\sat \s$ and $\depth(\uK')<\depth(\uK)$.  By construction of \proc, there
exists an instantiation of $\calR\in \{\to L, e\to L, \to R, e\to R, \K R\}$ with $\s$
whose result is $\s_1,\dots,\s_r$ ($r=2$ if $\calR\in \{\to L,e\to L\}$, $r=3$ otherwise)
and $\proc(\s_r)$ is a proof. Thus there exists $\s'\in\{\s_1,\dots,\s_{r-1}\}$ such that
$\rho'\sat\s'$. But by induction hypothesis the depth of $\proc(\s')$ is the minimum among
the models for $\s'$. Thus no such
$\uK'$ can exist.\medskip\\
%
{\bf $\uK$ is the result of performing Step~(\ref{proc:comparison}).} If $\proc(\s)$ returns
$\uM$, then the proof goes as in the case of Point~\ref{proc:glue-model}, otherwise the
proof goes as in the case of Point~\ref{proc:inv-model}.
\medskip\\
{\bf $\uK$ is the result of performing Step~\ref{proc:KL-return-invertible}.}
Proved by contradiction using the facts: (a) there
is no Kripke model satisfying $\s_2$; (b)  by induction hypothesis, $\uK$ is model of
minimum depth satisfying $\s_1$. \medskip\\
{\bf $\uK$ is the result of performing Step~\ref{proc:KL-return-glued}.}  By induction
hypothesis $\calM$ is minimal among the models satisfying $\s_2$ and if we are at this
step of \proc, then there is no model satisfying $\s_1$. The proof follows by
contradiction.\medskip\\
{\bf $\uK$ is the result of performing~\ref{proc:KL-return-comparison}.} In the previous
cases we have already proved that $U_1$ and $\calM$ satisfy the statement of the theorem.
\qed


\begin{example}
  \rm
  In previous section we have argued that the classical reflection principle \mbox{$\K A\to A$} is
  not provable in \Liel by using a combinatorial argument. Now, by using \proc and the
  completeness theorem, we formalise that \mbox{$\K A\to A$} is not a formula of the logic by
  providing a (counter)model that does not force it. We describe the steps performed by
  \mbox{$\proc(\seq{}{}{\mbox{$\K A\to A$}})$.}

  The sequent $\seq{}{}{\K A\to A}$, is used as actual parameter of \proc.
  Step~(\ref{proc:non-invertile-step:loop:calls}) is reached and the recursive calls
  $U_1=\proc(\seq{}{\K A}{A})$ and $U_2=\proc(\seq{}{\K A}{A})$ are performed. With
  $\seq{}{\K A}{A}$ as actual parameter, Step~\ref{proc:rule-K-L} is reached and the
  recursive calls $U_1=\proc(\eseq{\bot}{A}{A})$ and $U_2=\proc(\eseq{\bot}{A}{\bot})$
  are performed:
  \begin{itemize}
  \item with actual parameter $\eseq{\bot}{A}{A}$, procedure \proc reaches Step~\ref{proc:assiomi},
    thus the proof $\eseq{\bot}{A}{A} \ruleID$ is returned. Summarising:
    $U_1=\proc(\eseq{\bot}{A}{A})= \eseq{\bot}{A}{A} \ruleID$;
  \item with actual parameter $\eseq{\bot}{A}{\bot}$, procedure \proc reaches
    Step~\ref{proc:no-rule-to-apply} because $\eseq{\bot}{A}{\bot}$ is a flat
    sequent and a Kripke model is returned. Summarising
    $U_2=\proc(\eseq{\bot}{A}{\bot})=\langle\{\rho\}, \rho, \{(\rho,\rho)\},
    \{(\rho,\rho)\}, \{(\rho, A)\}\rangle$. 
  \end{itemize}
  After the two recursive calls at Step~\ref{proc:rule-K-L} terminate, since $U_1$ is a
  proof, Step~\ref{proc:rule-K-L:build-structure} is performed and from the Kripke model
  \[
    U_2=\langle \{\rho'\}, \rho', \{(\rho',\rho')\}, \{(\rho',\rho')\}, \{(\rho',A)
    \rangle
  \]
  the Kripke model $\uM=\kripke$ is defined as follows:
  \[
    S=\{\rho,\rho'\},\ \leq=\{(\rho,\rho),
    (\rho,\rho'),(\rho',\rho')\},\ \E=\{(\rho,\rho'),\ (\rho',\rho')\},\ \forza=\{(\rho',A)\}
  \]
  Since $U_1$ is a proof, $\calM$ is returned in
  Step~\ref{proc:rule-K-L:return-structure}.  At this point we have that the two calls
  $U_1=\proc(\seq{}{\K A}{A})$ and $U_2=\proc(\seq{}{\K A}{A})$ in
  Step~(\ref{proc:non-invertile-step:loop:calls}) return in both cases $\calM$. The sets
  $\inv$ and $\noninv$ are respectively updated at
  Steps~(\ref{proc:non-invertile-step:loop:first-update-U1})
  and~(\ref{proc:non-invertile-step:loop:first-update-U2}). Since $\noninv$ is not empty,
  Step~(\ref{proc:glue-step}) is performed and $\calM$ of depth three is built. Since
  $\inv=\{U_1\}$ and $U_1$ has depth two, we get that $\proc(\seq{}{}{\K A\to A})$ is the
  structure $U_1=\calM$ of depth two returned performing Step~(\ref{proc:comparison}).
\end{example}

\begin{example}
  \rm $\K A$ is invalid in $\iel$. We perform the call $\proc(\seq{}{}{\K A})$. Since the
  actual parameter is $\proc(\seq{}{}{\K A})$,
  Step~(\ref{proc:non-invertile-step:loop:calls}) is reached and the recursive calls
  $U_1=\proc(\eseq{}{}{A})$ and $U_2=\proc(\seq{}{}{A})$ are performed. Since
  $\eseq{}{}{A}$ and $\seq{}{}{A}$ are flat sequents, in both calls
  Step~\ref{proc:no-rule-to-apply} is reached and a Kripke model of depth one is
  returned. The sets $\inv$ and $\noninv$ are respectively updated at
  Steps~(\ref{proc:non-invertile-step:loop:first-update-U1})
  and~(\ref{proc:non-invertile-step:loop:first-update-U2}). Since $\noninv$ is not empty,
  Step~(\ref{proc:glue-step}) is performed and $\calM$ of depth two is built. Finally
  Step~(\ref{proc:comparison}) is reached and $U_1$ of depth one is returned as result of
  $\proc(\seq{}{}{\K A})$. Note that the model proving the invalidity in \iel of $\K A $
  has a single world that \E-reaches himself.
\end{example}

\begin{example}
  \rm
  $\K (A\lor B)\to (\K A \lor \K B)$ is invalid. There are two possible
  completed trees of sequents. A part of them is provided in the following (in the $R \to$
  application we have two equal sequents, thus we  show one).  
  
  {
    \newcommand{\axaa}{\seq{}{A}{A}~\ruleID}
    \newcommand{\axab}{\seq{}{A}{B}}
    \newcommand{\axba}{\seq{}{B}{A}}
    \newcommand{\axbb}{\seq{}{B}{B}~\ruleID}
    \newcommand{\deda}{ \regola{\seq{}{A\lor B}{A} }{\axaa & \axba}{\lor L}}
    \newcommand{\dedsx}{\regola{\eseq{}{A\lor B}{A, \K B}}{\vdots}{\K R}}
    \newcommand{\dedb}{ \regola{ \seq{}{\K(A\lor B)}{\K A, \K B} }
      { \dedsx &\deda}{\K R}}
    \newcommand{\dedc}{ \regola{ \seq{}{\K(A\lor B)}{\K A \lor \K B}}{\dedb}{ \lor R} }
    \newcommand{\dedd}{ \regola{ \seq{}{}{\K(A\lor B) \to  (\K A \lor \K B)}}{\dedc}{ \to R} }
    \[
      \dedd
    \]
  }
  {
    \newcommand{\axaa}{\seq{}{A}{A}~\ruleID}
    \newcommand{\axab}{\seq{}{A}{B}}
    \newcommand{\axba}{\seq{}{B}{A}}
    \newcommand{\axbb}{\seq{}{B}{B}~\ruleID}
    \newcommand{\deda}{ \regola{\seq{}{A\lor B}{B} }{\axab & \axbb}{\lor L}}
    \newcommand{\dedsx}{\regola{\eseq{}{A\lor B}{B, \K A}}{\vdots}{\K R}}
    \newcommand{\dedb}{ \regola{ \seq{}{\K(A\lor B)}{\K A, \K B} }{\dedsx & \deda}{\K R}}
    \newcommand{\dedc}{ \regola{ \seq{}{\K(A\lor B)}{\K A \lor \K B}}{\dedb}{ \lor R} }
    \newcommand{\dedd}{ \regola{ \seq{}{}{\K(A\lor B)\to (\K A \lor \K B)}}{\dedc}{ \to R} }
    \[
      \dedd
    \]
  }

  Note that there are two possible instantiations of $\K R$, this explains the two
   trees of sequents. In both cases, the tree of sequents ends with a
  flat sequent. Thus we get two models. The collection of the trees of sequents is a way
  to represent the recursive calls of \proc. The rightmost branch of the first tree
  of sequents ends with the flat sequent ${\seq{}{B}{A}}$. From the construction of \proc there is
  a model with one world that satisfies ${\seq{}{B}{A}}$. The same model satisfies $\seq{}{A\lor
    B}{A}$. Similarly for the second tree of sequents starting from the flat sequent
  $\seq{}{A}{B}$. Sequent  $\seq{}{\K(A\lor B)}{\K A, \K B}$ instantiates $\K R$ in two
  different ways, thus this is a backtracking point. The two models we have at hand are
  ``glued'' together using a new world in Step~\ref{proc:glue-step} and a model with three
  worlds is built. This model satisfies the remaining sequents in the two branches we have
  considered in the two trees of sequents.

\end{example}


\section{Refutational calculus for \iel}
\label{sec:riel}
\newcommand{\prema}{\{\seq{B}{\T,\G_1}{A}\}_{A\to B\in \G}}
\newcommand{\premb}{\{\seq{}{ A, \T,\G}{B}\}_{A\to B\in\D}}
\newcommand{\premc}{\{\seq{}{ \T,\G_2}{B}\}_{\K B\in\D}}
Usually we are interested in designing calculi that prove the validity of a given formula
$A$. As a result, in a model theoretic approach, the validity of $A$ is witnessed by a
proof and the invalidity is witnessed by a model that does not satisfy $A$.

This asymmetry between the validity, where a tree of sequents is returned and the
invalidity, where a relational structure is returned, can be adjusted by designing a
logical calculus to prove the invalidity.

To keep the two aspects aparted, calculi aimed to prove the invalidity are called {\em
  refutational calculi}. The proofs built with refutational calculi are called {\em
  refutations}. A formula provable in a refutational calculus is called {\em refutable}.
Since a refutable calculus aims to prove the invalidity, if a formula is refutable, then
there exists a world of a Kripke model that does not force it. This is the correctness of the
refutational calculus. We also want that if a formula $A$ is invalid, that is there exists
a world of a Kripke model that does not force $A$, then $A$ is refutable, that is provable in the
refutational calculus. This is the completeness of the refutational calculus.

The aim of this section is to present a refutational calculus for \iel. The work developed
to prove the correctness and the completeness of \Liel is useful also to prove correctness
and completeness of the refutational calculus \Riel provided in
Figure~\ref{fig:single-conclusion-rules-for-refutational-calculus}.
\begin{figure}[t]
  \footnotesize
  \[
    \renewcommand{\arraystretch}{1}
    \begin{array}{c}
      {\bf Axioms} \medskip\\
      \seq{\T}{\G}{\D}~~\ruleSAT\hspace*{4em}\eseq{\T}{\G}{\D}~~\ruleeSAT
      \medskip\\
      \textrm{ provided } \G\subseteq\PV,\,\D\subseteq\AT\textrm{ and } \G\cap\D=\emptyset %
      \medskip\\
      {\bf Rules}\smallskip\\
      \textrm{ (Proviso: rules apply  iff $\bot$ does not occur in the second compartment
      and}\\
      \textrm{the second and third compartments do not share  formulas)} \medskip\\
      \infer[(\land L)]{\seq{\T}{A \land B,\G}{\D}}
      {\seq{\T}{A, B, \G}{\D}} %
      \hspace{4em}
      \infer[(\land R_i)]{\seq{\T}{\G}{A_1 \land A_2, \D}}
      {\seq{\T}{\G}{A_i,\D} } %
      \medskip\\
      \infer[(\lor L_i)]{\seq{\T}{A_1 \lor A_2,\G}{\D}}
      {\seq{\T}{A_i,\G}{\D}} %
      \hspace{4em}
      \infer[(\lor R)]{\seq{\T}{\G}{A \lor B, \D}}
      {\seq{\T}{\G}{A,B,\D}} %
      \medskip\\
      \infer[(\to L_1)]{\seq{\T}{A \to B,\G}{\D}}
      {\seq{\T}{B,\G}{\D}}
      \hspace*{4em}
      \infer[(\to L_2)]{\seq{\T}{A \to B,\G}{\D}}
      {\seq{B,\T}{\G}{A,\D} }
      \medskip\\
      \regola{\seq{\T}{\K A,\G}{\D}}
      {\eseq{\false}{A,\G}{\D}}{\K L_1}%
                                  \hspace{4em}
      \regola{\seq{\T}{\K A,\G}{\D}}
      { \eseq{\false}{A, \T,\G}{\false}}{\K L_2}\raise 1.3ex\hbox{, provided $\D\subseteq\AT$} %
      \medskip\\
      \regola{\seq{\T}{\G}{A \to B,\D}}
      {\seq{\T}{A, \G}{B,\D} }{\to R_1}\hspace{4em}
      \regola{\seq{\T}{\K A_1,\dots,\K A_n,\G }{\K B,\D}}
      {\eseq{\T}{A_1,\dots, A_n,\G}{B,\D}}
      {\K R_1}\smallskip\\
      \textrm{\hspace*{15em}where $n\geq 0$}
      \medskip\\
      {
      \def\concl{\eseq{\T}{\K A, \G}{\D}}
      \def\prem{\eseq{\T}{A,\G}{\D}}
      \regola{\concl}{\prem}{e\K L}
      }
      \hspace*{4em}
      {
      \def\concl{\eseq{\T}{\G}{\K B,\D}}
      \def\prem{\eseq{\T}{\G}{B,\D}}
      \regola{\concl}{\prem}{e\K R}
      }
      \medskip\\
      \infer[(e\land L)]{\eseq{\T}{A \land B,\G}{\D}}
      {\eseq{\T}{A, B, \G}{\D}} %
      \hspace{4em}
      \infer[(e\land R_i)]{\eseq{\T}{\G}{A_1 \land A_2, \D}}
      {\eseq{\T}{\G}{A_i,\D} } %
      \medskip\\
      \infer[(e\lor L_i)]{\eseq{\T}{A_1 \lor A_2,\G}{\D}}
      {\eseq{\T}{A_i,\G}{\D} } %
                            \hspace{4em}
                            \infer[(e\lor R)]{\eseq{\T}{\G}{A \lor B, \D}}
                            {\eseq{\T}{\G}{A,B,\D}} %
      \medskip\\
      \infer[(e\to L_1)]{\eseq{\T}{A \to B,\G}{\D}}
      {\eseq{\T}{B,\G}{\D} }
      \hspace{4em}
      \infer[(e\to L_2)]{\eseq{\T}{A \to B,\G}{\D}}
      { \eseq{B,\T}{\G}{A,\D} }
      \medskip\\
      \regola{\eseq{\T}{\G}{A \to B,\D}}
      {\eseq{\T}{A, \G}{B,\D} }{e\to R_1}
      \medskip\\
      \regola{\seq{\T}{\G}{\D}}{\prema & \premb & \premc}{Glue}\smallskip\\
      \textrm{where }
      \G\subseteq\PV\cup\{A\to B\in\calL\}\cup\{\K A\in\calL\},
      \D\subseteq\AT\cup\{A\to B\in\calL\}\cup\{\K A\in\calL\},\smallskip\\
      \G_1= \G\setminus\{A\to B\}\textrm{ and }
      \G_2=(\G\setminus\{\K A\in \G\})\cup \{A | \K A\in \G\}
      \medskip\\
      \regola{\eseq{\T}{\G}{\D}}{\prema & \premb}{eGlue}\smallskip\smallskip\\
      \textrm{where }
      \G\subseteq\PV\cup\{A\to B\in\calL\},
      \D\subseteq\AT\cup\{A\to B\in\calL\} \textrm{ and }
      \G_1= \G\setminus\{A\to B\}\\
    \end{array}
  \]
  \caption{The refutational calculus \Riel for \iel.}
  \label{fig:single-conclusion-rules-for-refutational-calculus}
\end{figure}                     
\clearpage
We remark that \Liel and \Riel share the same object language. The names of the rules of
\Riel are derived from those of \Liel.  The ideas behind the rules of \Riel are the fact
that, if a given formula $A$ is not valid, then by them we construct a model in which a
world does not force $A$. In practice we are exploiting
Proposition~\ref{prop:satisfiability}. We remark that ($\land L$), ($\lor R$), ($e\K L$),
($e\K R$), ($e\land L$) and ($e\lor R$) are the same rules of \Liel. This is related to the fact
that a model $\uK$ satisfies the premise iff $\uK$ satisfies the conclusion, or,
equivalently, the conclusion is valid iff the premise is valid.  \Liel
handles the disjunction of the left by means of the two-premise rule ($\lor L$). \Riel has
the single-premise rules ($\lor L_1$) and ($\lor L_2$). Thus disjunctions on the left are
backtracking points in \Riel proof search. The motivation for this is clear at the light
of Proposition~\ref{prop:satisfiability}: roughly speaking, to prove the invalidity of a
disjunction on the left is sufficient to prove that one of the disjuncts is invalid.  The
same applies to conjunctions on the right. As regards implications and $\K$-formulas they
are handled both by the {\em ad-hoc} rules ($\to R_1$), ($\K R_1$), ($e\to L_1$), ($e\to L_2$),
($e\to R_1$) and by the {\em collective} rules ($Glue$) and ($eGlue$). This can be informally
explained by the rules of \Liel and by Proposition~\ref{prop:satisfiability}. As an
example, \Riel handles implications on the left of sequents of the kind $\seq{\T}{\G}{\D}$
by the rules ($\to L_1$), ($\to L_2$) and ($Glue$). This is related to the fact that \Liel
handles the implications on the left by a rule with three premises and, by
Proposition~\ref{prop:satisfiability}, we know that if the first or second premise of
$\to L$ is satisfied, then the conclusion of ($\to L$) is satisfied. In general we cannot
draw the same conclusion for the third premise of ($\to L$), and thus \Riel needs ($Glue$). Same
remarks apply for the other cases.  Clearly, the duality between \Liel and \Riel is
related to the fact that \Liel is designed to prove validity, \Riel to prove invalidity.

%
To prove that \Riel is correct, we have to show that if a formula $A$ is \Riel-provable,
then $A$ is not valid in \iel. First of all we have to fix the definitions given in
Section~\ref{sec:definitions} for \Liel to the case of \Riel. We remark that in \Riel the notions
of axiom and flat are switched w.r.t. \Liel. 
In \Riel the axioms are $\ruleSAT$ and $\ruleeSAT$ and the flat sequents contains $\bot$ in
the second compartment or their second and third compartment share formulas. 

Now we are ready to prove the correctness of \Riel.  This requires to prove that from a
\Riel-proof $\pi$ of $\s$ we can define a Kripke model for \iel. After that, by
Proposition~\ref{prop:riel:satisfiability} we prove that the root of the Kripke model
extracted from $\pi$ satisfies $\s$. Both these tasks have been already done:

\begin{proposition}[Satisfiability]
  \label{prop:riel:satisfiability}
  Let $s$ be a sequent in the object language of \Riel. Let $\pi$ be a \Riel proof
  of $\s$. Then there exists a Kripke model $\uK=\kripke$, such that $\rho\sat\s$.
\end{proposition}
{\em Proof:} we proceed by induction on the depth of $\pi$. Note that this proposition
corresponds to prove Theorem~\ref{th:kripke-model} and
Proposition~\ref{prop:satisfiability} for \Liel, where $\pi$ has the role of \proc. Thus
we exploit the work already done and we provide only a proof sketch.
\smallskip\\
{\em Base:} no rule is applied. Thus $\pi$ coincides with $\s$. Sequent $\s$ is an axiom
of \Riel. This means that $\s$ is flat sequent of \Liel. Let $\uK$ be defined as in
Step~\ref{no-rule-to-apply} of \proc. By Theorem~\ref{th:kripke-model}, $\uK$ is a Kripke
model. Now the statement of the theorem is proved as in the base case of
Proposition~\ref{prop:satisfiability}.
\smallskip\\
{\em Induction:} we proceed by cases according the rule $\rmR$ that $\s$ instantiates in
the construction of $\pi$. We assume that for every sequent $\s'$ in the set resulting
from the instantiation of $\rmR$ with $\s$, there exists a Kripke model $\uK'$ with root
$\rho'$ such that $\rho'\sat\s'$.  We notice that for every possible value of $\rmR$ we
have already proved the result in Theorem~\ref{th:kripke-model} and
Proposition~\ref{prop:satisfiability}.
\smallskip\\
As regards the rules in the sets $\{(\land L), (\lor R), (e\land L), (e\lor R), (e\K L), (e\K R) \}$,
$\{(\lor L_i), (e\lor L_i), (\land R_i), (e\land R_i)\}$,
$\{(\to L_1),(\to L_2), (\to R), (e\to L_1), (e\to L_2), (\K R)\}$ and $\{\K L_1\}$ we extract a
structure $\uK$ that coincides with $\uK'$ respectively following the 
Steps~(\ref{proc:alpha-step}), \ref{proc:beta-step:exactly-one}, (\ref{proc:inv-model})
and \ref{proc:KL-return-invertible} of \proc.  Theorem~\ref{th:kripke-model} proves that
$\uK$ is a model immediately.
As regards the satisfiability, in Proposition~\ref{prop:satisfiability}
see the proof respectively for
{\bf $\mathbfcal{K}$ is the result of performing Step~(\ref{proc:alpha-step})},
{\bf $\mathbfcal{K}$ is the result of performing Step~\ref{proc:beta-step:exactly-one}},
{\bf $\mathbfcal{K}$ is the result of performing Step~(\ref{proc:inv-model})}
and
{\bf $\mathbfcal{K}$ is the result of performing Step~\ref{proc:KL-return-invertible}.}
\smallskip\\
As regards the remaining rules:\smallskip\\
%
%
%
Rule $\K L_2$: from $\pi$ the structure $\uK$ is defined as in
Step~\ref{proc:KL-return-glued} of \proc using $\uK'$. Structure $\uK$ is proved to fulfil
the definition of Kripke model as in Theorem~\ref{th:kripke-model}.  As regards the
satisfiability, we proceed as in the proof of Proposition~\ref{prop:satisfiability} for
the case {\bf $\mathbfcal{K}$ is the result of performing
  Step~\ref{proc:KL-return-glued}.}\smallskip\\
Rules $(Glue)$ and $(eGlue)$: the case $\s=\seq{\T}{\G}{\D}$ and
$\s=\eseq{\T}{\G}{\D}$ are analogous. The models obtained by induction from the premises of the
rules are glued together as in Step~(\ref{proc:glue-model}) of \proc to define
$\uK$. Now, we proceed as in Theorem~\ref{th:kripke-model} to prove that $\uK$ is a Kripke
model and as in {\bf $\mathbfcal{K}$ is the result of performing
  Step~(\ref{proc:glue-model})} of Proposition~\ref{prop:satisfiability} to prove the
satisfiability.
\qed
\begin{theorem}[Correctness of \Riel]
  \label{th:riel-correctness}
  Let $A\in\calL$. If \Riel proves $\seq{}{}{A}$,  then $A$ is not valid in \iel. 
\end{theorem}
{\em Proof:} By Proposition~\ref{prop:riel:satisfiability}, from a proof of $\Riel$ we can
extract a Kripke model $\uK=\kripke$ such that $\rho\sat \seq{}{}{A}$. By the meaning of
satisfiability of a sequent follows that $\rho\nonforza A$ and hence $A$ is not valid in
\iel. \qed

As regards the completeness of \Riel, we must show that every unprovable (irrefutable)
formula in \Riel is valid in \iel. To this aim we exploit the work already done for \Liel
providing the procedure \priel that given $\seq{}{}{A}$, returns a proof of \Liel if $A$
is valid and a proof of \Riel if  $A$ is not valid. Procedure \priel is a
rewriting of procedure \proc, where for terseness we disregard the part related to
minimality.


\noindent
{\sc Procedure $\priel(\s)$} 
\begin{enumerate}[label=(\arabic*),ref=\arabic*]

\item \label{priel:assiomi}
  if $\s$ can instantiate a \Liel axiom, then return the \Liel proof \[\s~(\rmR)\] where
  $\rmR\in\{\mathbf{Irr}, \mathbf{Id},\mathbf{eIrr},\mathbf{eId}\}$;
  
\item\label{priel:no-rule-to-apply}
  if $\s$ can instantiate a \Riel axiom, then return the \Riel proof \[\s~(\rmR)\] where
  $\rmR\in\{\mathbf{Sat},\mathbf{eSat}\}$;
  
\item \label{priel:alpha-step} if $\s$ can instantiate
  $\rmR\in\{\land L, \lor R,e\land L, e\lor R, e\K L, e\K R \}$, then let $\s_1$ be the
  result of an instantiation of $\rmR$ with $\s$. Let $U=\priel(\s_1)$. Return
  $\apply{U}{\s}{(\rmR)}$, where $\apply{U}{\s}{(\rmR)}$ is a \Liel or \Riel proof
  according to the kind of $U$;
  
\item \label{priel:beta-step} if $\s$ can instantiate $\rmR\in\{\lor L, \land R, e\lor L,
  e\land R \}$, then let
  $\s_1$ and $\s_2$ be the result of an instantiation of $\rmR$ with $\s$.
  For $i=1,2$, let $U_i=\priel(\s_i)$;
  \begin{enumerate}
  \item If $U_1$ and $U_2$ are \Liel proofs, then return $\regola{\s}{U_1 & U_2}{\rmR}$;
  \item \label{priel:beta-step:first} if $U_1$ is a \Riel proof, then return
    $\regola{\s}{U_1 }{\rmR_1}$;
  \item \label{priel:beta-step:second} return
    $\regola{\s}{U_2}{\rmR_2}$;
  \end{enumerate}
\item \label{priel:non-invertile-step} if $\s$ can instantiate one of the rules in
  $\{(\to R)$, ($\to L$), ($\K R$), ($e\to R$), ($e\to L)\}$ then:
  \begin{enumerate}[label=(\alph*),ref=\theenumi.\alph*]
  \item let $\noninv=\emptyset$;
  \item
    \label{priel:non-invertile-step:loop}
    for every $\rmR\in \{\to R, \K R, e\to R\}$, for every possible instantiation of
    $\rmR$ with $\s$:
    \begin{enumerate}[label=(\roman*),ref=\theenumii.\roman*]
    \item \label{priel:non-invertile-step:loop:calls}
      let $\s_1$ and $\s_2$ the result of the instatiation of $\rmR$ with $\s$, let
      $U_1=\priel(\s_1)$ and $U_2=\priel(\s_2)$;
    \item \label{priel:proof-step} if $U_1$ and $U_2$ are \Liel proofs, then return
       $\regola{\s}{U_1 & U_2}{\rmR}$;
    \item \label{prile:non-invertile-step:loop:first-update-U1}
      if $U_1$ is a \Riel proof, then return $\regola{\s}{U_1}{\rmR_1}$;
    \item
      \label{priel:non-invertile-step:loop:first-update-U2}
      Collect the \Riel proof $U_2$: $\noninv=\noninv\cup\{U_2\}$
    \end{enumerate}
  \item for every $\rmR\in \{\to L, e\to L\}$, for every possible instantiation of
    $\rmR$ with $\s$:    
    \begin{enumerate}[label=(\roman*),ref=(\roman*)]
    \item let $\s_1$, $\s_2$ and $\s_3$ the result of the instantiation and let 
      $U_1=\priel(\s_1)$, $U_2=\priel(\s_2)$ and $U_3=\priel(\s_3)$;
    \item if $U_1$, $U_2$ and $U_3$ are \Liel proofs, then return
      $\regola{\s}{U_1 & U_2 &U_3}{\to L}$;
    \item if $U_1$ is a \Riel proof, then return $\regola{\s}{U_1}{\to L_1}$;
    \item if $U_2$ is a \Riel proof, then return $\regola{\s}{U_2}{\to L_2}$;
    \item collect the \Riel proof $U_3$: $\noninv=\noninv\cup\{U_3\}$;
      
    \end{enumerate}

  \item \label{priel:glue-step} Let $\s_1,\dots,\s_n$ be an enumeration of the \Riel proofs
    in $\noninv$. Return the \Riel proof $\regola{\s}{\s_1\dots\s_n}{\rmR}$, where
    $\rmR\in\{(Glue), (eGlue)\}$;
          
  \end{enumerate}
  
\item
  \label{priel:rule-K-L}
  let $\s_1$ and $\s_2$ be the result of an
  instantiation of $\K L$ with $\s$. Let $U_i=\priel(\s_i)$, for $i=1,2$;
  \begin{enumerate}
  \item If $U_1$ and $U_2$ are \Liel proofs, then return $\regola{\s}{U_1& U_2}{\K L}$;
  \item \label{priel:KL-return-invertible}%
    if $U_1$ is a \Riel proof, then return $\regola{\s}{U_1}{\K L_1}$;
  \item
    \label{priel:rule-K-L:build-structure}
    return \Riel proof $\regola{\s}{U_2}{\K L_2}$;
  \end{enumerate}
\end{enumerate}
{\sc End Procedure \priel}\medskip\\
%
By induction on the number of recursive calls, it is easy to show that given a sequent
$\s$, \priel returns a \Liel proof of $\s$ or a \Riel proof of $\s$.  By the correctness
of \Liel and \Riel we can get in a row the completeness of \Riel, which is the result that
we have to prove, but also, in another form the completeness of \Liel:
\begin{theorem}[Completeness of \Riel (and \Liel)]
  Let $A$ be a formula. If $\priel(\seq{}{}{A})$ returns a \Riel proof, then $A$ is not
  valid, otherwise $A$ is valid
\end{theorem}
{\em Proof:} we have already remarked that given a sequent, \priel returns a
\Liel proof or a \Riel proof. Let $A$ be an invalid formula. By the correctness of \Liel,
\priel cannot return a proof of \Liel for $\seq{}{}{A}$, thus \priel returns a proof of
\Riel. Thus we have proved that all the invalid formulas have a proof in the refutational
calculus \Riel. \qed Note that with a dual argument we can prove the completeness of
\Liel, by using \priel and the correctness of \Riel.

As regards procedure \priel, we remark that for every possible instantiation of $\s$
with $\{(\to R)$, ($\to L$), ($\K R$), ($e\to R$) ($e\to L)\}$, we consider the \Riel proof
$\priel(\s_r)$. Note that the collection of the sequents $\s_r$ coincides with the result of
instantiating $(Glue)$ or $(eGlue)$ with $\s$.  \Riel proofs collected in $\noninv$ are
the result of $\priel(\s_r)$.

\begin{example}\rm
  \Riel proves the invalidity of classical reflection $\K A\to A$. 

    \newcommand{\seqa}{\eseq{}{A}{}~\ruleeSAT}
    \newcommand{\axa}{\eseq{}{A}{A}\ruleID}
    \newcommand{\seqb}{\seq{}{\K A}{A}}
    \newcommand{\seqc}{\seq{}{}{\K A\to A}}
    \newcommand{\deda}{\regola{\seqb}{ \seqa}{\K L_2}}
    \[
      \regola{\seqc}{\deda}{\to R_1}
    \]

  Note that in this proof we have two backtracking points. As a matter of fact, seqnet $\seqc$
  can instantiate rules ($\to R_1$) and ($\to R_2$) and sequent $\seqb$ can instantiate
  rules $(\K L_1)$ and $(\K L_2)$. 
\end{example}

\begin{example}
  \rm \Riel proves the invalidity of $\K (A\lor B)\to (\K A \lor \K B)$. There are two
  possible completed trees of sequents that are provided in the following (where in the
  ($R \to$) application we have two equal sequents, thus we show one branch only).
  
  {
    \newcommand{\axaa}{\seq{}{A}{A}~\ruleID}
    \newcommand{\axab}{\seq{}{A}{B}~\ruleSAT}
    \newcommand{\axba}{\seq{}{B}{A}~\ruleSAT}
    \newcommand{\axbb}{\seq{}{B}{B}~\ruleID}
    \newcommand{\deda}{ \regola{\seq{}{A\lor B}{A} }{\axba}{\lor L_2}}
    \newcommand{\dedsx}{\regola{\seq{}{A\lor B}{B}}{\axab}{\lor L_1}}
    \newcommand{\dedb}{ \regola{ \seq{}{\K(A\lor B)}{\K A, \K B} }
      { \dedsx &\deda}{Glue}}
    \newcommand{\dedc}{ \regola{ \seq{}{\K(A\lor B)}{\K A \lor \K B}}{\dedb}{ \lor R} }
    \newcommand{\dedd}{ \regola{ \seq{}{}{\K(A\lor B) \to  (\K A \lor \K B)}}{\dedc}{ \to R_1} }
    \[
      \dedd
    \]
  } We remark that in proof construction there are some backtracking points:
  \begin{enumerate}
  \item $\seq{}{}{\K(A\lor B) \to  (\K A \lor \K B)}$ can also instantiate rule
    $(Glue)$. If one is not interested in minimality, this choice makes no difference;
  \item $\seq{}{\K(A\lor B)}{\K A \lor \K B}$ can also instantiate $(\K L_1)$ ;
  \item $\seq{}{\K(A\lor B)}{\K A, \K B}$ can also instantiate $(\K L_1)$, and, in two
    ways, $(\K R_1)$;
  \item $\seq{}{A\lor B}{A}$
    can also instantiate $(\lor L_1)$;
  \item $\seq{}{A\lor B}{B}$ can also instantiate $(\lor L_2)$. 
  \end{enumerate}
\end{example}
By the work developed in Section~\ref{sec:liel:completeness} and in
Proposition~\ref{prop:riel:satisfiability} it should be straightforward how-to extract a
Kripke model whose root does not satisfy a sequent provable in \Riel.

We conclude by remarking that by inspection of the rules, the depth of all
$\Riel$-trees $\calT$ is bounded by the number of connectives occurring in the sequent 
in the root  of $\calT$.


\section{Calculus for \ielmeno}
\label{sec:lielmeno}
Logic $\iel^-$ lacks of axiom~\ref{Ax4}. Semantically this means that \mbox{property~\ref{Im3}}
on Kripke models does not hold.

We show that we get the complete calculus $\Liel^-$ for \ielmeno by removing rule $\K L$ from the
logical apparatus for \Liel. The intuition is related to \proc, that uses $\K L$ in
Step~\ref{proc:rule-K-L} when no other rule is applicable. When \proc performs
Step~\ref{proc:rule-K-L}, the following statement ($\calS$) on $\seq{\T}{\G}{\D}$ holds:
\begin{enumerate}
\item[($\calS$)]$\D\subseteq \AT$, the second compartment is included in $\{\K A\in \cal L\}\cup \PV$ and
the intersection between second and third compartment is empty.  
\end{enumerate}
The following formally justifies our choice and is part of the completeness theorem:
\begin{proposition}
  \label{prop:ielmeno:flat}
  The sequent $\s=\seq{\T}{\K A,\G}{\D}$ that satisfies statement~($\calS$) is
  satisfiable by the root of the Kripke model $\uK=\kripke$, where $S=\{\rho\}$,
  $\leq=\{(\rho,\rho)\}$, $\E=\emptyset$, $\forza=\{(\rho,p) | p\in \G\cap\PV\}$.   
\end{proposition}
{Proof:}
Since $\G$
and $\D$ do not share any element and by definition of forcing, $\rho\forza\bigvee\D$. By
definition of forcing, for every $p\in\PV$, if $p\in\G$, then $\a\forza p$. By the fact
that $\E=\emptyset$ and the semantics of $\K$, for every $\K A\in \G$, $\a\forza \K
A$. Since $\uK$ has a single world it is immediate that $\rho$ satisfies the first
compartment. Hence we have proved that $\rho$ satisfies $\s$.\qed

Since there is no rule to handle directly $\K$-formulas on the left, we have a new
definition of {\em flat} sequent for calculus $\Liel^-$ that extends the definition of
flat sequent for \Liel. Flat sequents for $\Liel^-$ are  $\E$-sequent
$\eseq{\T}{\G}{\D}$ fulfilling $\G\subseteq \PV$, $\D\subseteq\AT$ and
$\G\cap\D=\emptyset$ or sequents $\seq{\T}{\G}{\D}$ fulfilling
$\G\subseteq \{\K A\in\calL\}\cup \PV$, $\D\subseteq\AT$ and $\G\cap\D=\emptyset$.
As usual, the definition of flat sequent is tailored to characterise a non-axiom sequent that does
not instantiate any rule of the calculus.

We use the new definition of flat sequent to provide procedure $\proc^-$ that returns a
proof of $\Liel^-$ or a Kripke model for \ielmeno.  Because of the close relationships
between \Liel and $\Liel^-$, we have that   Procedure $\proc^-$ is a
slight modification of $\proc$:
\begin{itemize}[label=-]
\item erase Step~\ref{proc:rule-K-L} from procedure \proc ;
\item replace Step~\ref{proc:no-rule-to-apply} of \proc with the following two steps:
  \begin{enumerate}[label=2(\alph*)]
  \item\label{procmeno:no-rule-to-apply-eseq}%
    if $\s$ is a flat sequent  $\eseq{\T}{\G}{\D}$, then return the
    structure
    \[
      \uK=\langle  S, \rho, \leq, \E, \forza\rangle
    \]
    where $S$ = $\{\rho\}$, $\leq$ = $\{(\rho,\rho)\}$, $\E=\{(\rho,\rho)\}$  and $\forza$ =
    $\{\rho\}\times (\G\cap\PV)$;

  \item \label{procmeno:no-rule-to-apply-seq}%
    if $\s$ is a flat sequent  $\seq{\T}{\G}{\D}$, then return the
    structure
    \[
      \uK=\langle  S, \rho, \leq, \E, \forza\rangle
    \]
    where $S$ = $\{\rho\}$, $\leq$ = $\{(\rho,\rho)\}$, $\E=\emptyset$  and $\forza$ =
    $\{\rho\}\times (\G\cap\PV)$;
  \end{enumerate}
  
\end{itemize}

Now, the proof that for every invalid formula $\proc^-$ returns a Kripke model for \ielmeno of minimal
depth follows from Proposition~\ref{prop:ielmeno:flat} and is analogous to the proofs
given for Theorem~\ref{th:kripke-model}, Theorem~\ref{th:completeness} and Theorem~\ref{th:minimality}.

\begin{example}
  \rm Calculus $\Liel^-$ does not prove~\ref{Ax4} $\K A\to \neg\neg A$.
  Rules~$\K L$ and $\to R$ are instantiated with $\D=\{\bot\}$ and $\T=\emptyset$. Thus
  the premises are equal and to save space we only show one branch.
  \newcommand{\axa}{\seq{}{\K A,\bot}{\bot}~\ruleIRR}
  \newcommand{\axb}{\seq{\bot}{\K A}{A,\bot}}
  \newcommand{\axc}{\seq{\bot}{\K A}{A}}
  \newcommand{\seqa}{\eseq{}{A,\neg A}{}}
  \newcommand{\deda}{\regola{\seqb}{\axa & \axb & \axc}{e\to L}}
  \newcommand{\seqb}{\seq{}{\K  A,\neg A}{}}
  \newcommand{\seqc}{\seq{}{\K A}{\neg\neg A}}
  \newcommand{\dedc}{\regola{\seqc}{\deda}{\to  R}}
  \newcommand{\seqd}{\seq{}{}{\K A\to\neg\neg A}}
  \newcommand{\dedd}{\regola{\seqd}{\dedc}{\to R}}
  \[
    \dedd
  \]

  From the flat sequent on the middle we extract the model of minimum depth that does not
  forces the given formula: it is a model with a single world that does not force any
  propositional variable and $\E=\emptyset$.
\end{example}

To conclude the section, we remark that to get a refutational calculus $\Riel^-$ for
\ielmeno we  modify the refutational calculus \Riel for \Liel as follows:
\begin{itemize}[label=-]
\item the sequents fulfilling the property stated in~($\calS$) are axioms. Thus a sequent
  of the kind
  \[
    \begin{array}{c}
      \seq{\T}{\G}{\D}~(\mathbf{kSat})
      \\
      {\rm provided\ } \G\subseteq \{\K A\in\calL\}\cup\PV,~\D\subseteq\AT~{\rm and}~\G\cap\D=\emptyset.
    \end{array}
  \]
  is an axiom  of the refutational calculus $\Riel^-$;
\item the rules $\K L_1$ and $\K L_2$ of $\Riel$ are not rules of $\Riel^-$.
\end{itemize}
Now, by proceeding as for the case of \Riel, we can prove that $\Riel^-$ is a calculus for
the invalidity in \ielmeno.

\section{Conclusions and Future Works}
\label{sec:conclusions}
In this paper we have presented sequent calculi to prove validity and invalidity for the
intuitionistic propositional logics of belief and knowledge \iel and $\iel^-$. As for the
case of propositional intuitionistic logic~\cite{FerFioFio:2013}, we have shown that \iel
and $\iel^-$ have terminating calculi whose trees all have depth bounded by the number of
connectives in the formula to be proved and obey the subformula property. For invalid
formulas, our calculi allow us to get Kripke models of minimal depth. Compared
with~\cite{FerFioFio:2013}, the particular Kripke semantics characterising \iel and
$\iel^-$ requires an extension of the object language employed by the logical apparatus.
As in the case of~\cite{FerFioFio:2013}, the sequents are not standard and they have some
features related to nested sequents~\cite{Fitting:2014}. Roughly specking, a single
sequent to be satisfied requires that the formulas in the first compartment are satisfied in
all the successors of the world the second and third compartment refer to.

A possible future investigation is to get terminating calculi for intuitionistic epistemic
logic following the ideas in~\cite{Vorobiev:58} and~\cite{Hudelmaier:92}. These calculi do
not have the subformula property but they use standard sequents. A feature of the calculi
in~\cite{Hudelmaier:92} is that to get proofs in linear depth it can be required to
introduce new propositional variables. Finally, a further line of work is to apply the results
in~\cite{FerFioFio:2015} where it is showed a terminating strategy for the sequent
calculus {\bf G3i}. The strategy builds finite trees by using the information in the sequent at hand only,
no history mechanisms are required. In this case the depth of the returned proofs is
quadratic in the number of connectives occurring in the formula to be proved.

\end{document}